\documentclass[preprint,3p,authoryear]{elsarticle}

\usepackage{tikz}
\usepackage{booktabs}
\usepackage{float}
\usepackage{amsmath}
\usepackage{amsthm}
\usepackage{amsmath,amssymb,amsfonts,mathrsfs,hyperref,microtype}
\usepackage{graphicx}
\usepackage{graphicx,mathabx}
\usepackage[font=small,skip=0pt]{caption}
\graphicspath{ {images/} }
\usepackage{enumitem}  
\usepackage{chngcntr}
\newtheorem{theorem}{Theorem}[]
\newtheorem{lemma}{Lemma}[]
\newtheorem{corollary}{Corollary}[]
\newtheorem{proposition}{Proposition}[]
\theoremstyle{definition}

\newtheorem{exmp}{Example}[]
\newtheorem*{exmp*}{Example}
\theoremstyle{remark}

\newcommand{\BS}{\boldsymbol}
\newcommand\eop{\hfill$\Box$\\}

\newcommand{\Rmnum}[1]{\expandafter\@slowromancap\romannumeral #1@}
\newcommand{\trp}{{\sf T}}

\journal{}
\makeatletter
\def\ps@pprintTitle{%
   \let\@oddhead\@empty
   \let\@evenhead\@empty
   \def\@oddfoot{\reset@font\hfil\thepage\hfil}
   \let\@evenfoot\@oddfoot
}

\begin{document}
\begin{frontmatter}
\author[]{ Osama Idais\corref{cor1}}
\cortext[cor1]{Corresponding author}
\ead{osama.idais@ovgu.de. }
\author[]{Rainer Schwabe}
\ead{ rainer.schwabe@ovgu.de.}
\address{\small Institute for Mathematical Stochastics,  Otto-von-Guericke University  Magdeburg,\\ \small Universit\"atsplatz 2, 39106 Magdeburg, Germany}
 \title{In- and Equivariance for Optimal Designs in Generalized Linear Models: The Gamma Model}

\begin{abstract}
We give an overview over the usefulness of the concept of equivariance and invariance in the design of experiments for generalized linear models. 
In contrast to linear models here pairs of transformations have to be considered which act simultaneously on the experimental settings and on the location parameters in the linear component.
Given the transformation of the experimental settings the parameter transformations are not unique and may be nonlinear to make further use of the model structure.
The general concepts and results are illustrated by models with gamma distributed response.
Locally optimal and maximin efficient design are obtained for the common \textit{D}- and \textit{IMSE}-criterion.
\end{abstract}

\begin{keyword}
optimal design \sep invariance \sep equivariance\sep  generalized linear model\sep \textit{D}-criterion\sep \textit{IMSE}-criterion \sep maximin efficiency\
\end{keyword}

\end{frontmatter}

\section{Introduction}
\label{sec-1}
Generalized linear models are a powerful tool to analyze data for which the standard linear model approach is not adequate. 
The idea of generalized linear models goes back to \citet{10.2307/2344614} and the their concept is comprehensively presented in the monograph by \citet{mccullagh1989generalized}.
The statistical analysis is well developed in generalized linear models and there is also a considerable amount of literature on optimal design in this situation (see e.\,g.\ \citet{atkinson2015designs}).

In generalized linear models the performance of a design depends not only on the experimental settings but, in contrast to linear models, also on the value of the underlying parameters.
Even worse, also the optimal solutions for designs depend on the parameters.
As those are commonly unknown at the design step these parameters have to be guessed (estimated) prior to the experiment which leads to local optimality (see \citet{chernoff1953}). 
This approach has been frequently employed (see \citet{10.2307/2346142}, \citet{ATKINSON1996437}, \citet{10.2307/24308852}, \citet{10.2307/24309568},  \citet{yang2009},\citet{tong2014}, \citet{GAFFKE2019}, besides others), and it provides at least a benchmark for the quality of a design.  

To overcome the parameter dependence robust criteria have been proposed which either impose a prior weight on the parameters (Bayesian design, see e.\,g.\ \citet{AtkinsonDonevTobias}, ch.~18) or choose a minimax approach over a parameter region of interest (maximin efficiency, see e.\,g.\ \citet{10.2307/2345917}, \citet{10.2307/2677110}, \citet{DETTE20064397}, \citet{grasshoff2008optimal}).

In any case, the construction of optimal designs for generalized models is difficult and often numerical algorithms are used to find solutions.
However, to reduce the complexity of the search, it is advisable to make use of underlying symmetries (``invariance'') in the design problem.
The concept of invariance is well established for design optimization in linear models (see \citet{pukelsheim1993optimal}, ch.~5, or \citet{schwabe1996optimum}. ch.~3), but has not been used much in generalized linear models.
This is mainly due to the fact that local optimality criteria lack symmetries, in general, because of their parameter dependence.
The underlying concept of equivariance, however, has been profitably employed in generalized linear models with its most famous representative being the canonical transformations in \citet{10.2307/2346142}.
For the concept of invariance also symmetries in the parameters are required which concur with the symmetries in the experimental settings to account for the parameter dependence.  
In particular, Bayesian or maximin efficiency criteria can show additional symmetries in their prior or in the parameter region of interest (see \citet{radloff2016invariance})

Based on the approach sketched in \citet{radloff2016invariance} we develop the use of equivariance and invariance in the optimal design for generalized linear models in the present paper step by step and illustrate each step by the particular case of the gamma model, where the responses are gamma distributed. 
Such models are appropriate for many real life data from psychology, ecology or medicine. 
For the gamma model an additional reduction is possible by scaling of the parameters which leads to more complicated symmetry structures.

The paper is organized as follows.
 In Section~\ref{sec-2}, we introduce the model assumptions and the design criteria. 
 In Section~\ref{sec-3}, we discuss the concept of equivariance under standard linear transformations of the parameters and show how optimal designs can be transferred from one experimental region to another. 
 In Section~\ref{sec-3-extended}, the concept of equivariance is extended to nonlinear transformations of the parameters.
 In Section~\ref{sec-4}, the general concept of invariance is introduced and optimal designs are obtained for various situations.  
 Finally, Section~\ref{sec-5} concludes the paper with a short discussion. 

\section{Basics: Model specification, information, and design} 
\label{sec-2}

We consider a response variable $Y$ for which the dependence on a (potentially multi-dimensional) covariate $\mathbf{x}$ can be described by a generalized linear model.
This means that the distribution of $Y$ comes from a given exponential family and the mean $\mu=\mathrm{E}(Y)$ is related to the linear component $\mathbf{f}(\mathbf{x})^\trp\BS{\beta} = \sum_{j=0}^{p-1}\beta_j f_j(\mathbf{x})$ by a one-to-one link function.
In the linear component  $\mathbf{f}(\mathbf{x})=(f_0(\mathbf{x}),\ldots,f_{p-1}(\mathbf{x}))^\trp$ is a $p$-dimensional vector of given regression functions $f_0(\mathbf{x}),\ldots,f_{p-1}(\mathbf{x})$ and $\BS{\beta}=(\beta_0,\ldots,\beta_{p-1})^\trp$ is a  $p$-dimensional vector of parameters $\beta_0,\ldots,\beta_{p-1}$ to be estimated.
Traditionally the link function maps the mean to the linear component (see \citet{mccullagh1989generalized}, ch.~2).
For analytical purposes, however, it is more convenient to describe the dependence of the mean on the linear component,
\begin{equation}
	\label{eq2-1}
	\mu=\mu(\mathbf{x};\BS{\beta})=\eta(\mathbf{f}(\mathbf{x})^\trp\BS{\beta}),
\end{equation}
where $\eta$ is the inverse of the link function.
For example, for the $\log$ link $\eta$ is the exponential function.

As a particular case and for illustrative purposes we consider the gamma model with gamma distributed responses $Y$ with density given by $f_Y(y)=y^{\kappa-1}\exp(-y/\theta)/(\theta^\kappa\Gamma(\kappa))$, 
where $\kappa>0$ denotes the shape parameter and $\theta>0$ is the scale parameter.
Then $Y$ has expectation $\mu=\kappa\theta$.
For the link function we assume the commonly used inverse link $\kappa/\mu=\mathbf{f}(\mathbf{x})^\trp\BS{\beta}$, where the shape parameter $\kappa$ is supposed to be a fixed nuisance parameter (see \citet{atkinson2015designs}).
For example, $\kappa=1$ gives the family of exponential distributions, or for fixed integer $\kappa$ one obtains a family of certain Erlang distributions.

The inverse link is equal to the canonical link $-\kappa/\mu=\mathbf{f}(\mathbf{x})^\trp\BS{\beta}$ up to the minus sign (see \citet{mccullagh1989generalized}, ch.~3).
Its inverse $\eta$ is given by
\begin{equation}
	\label{eq2-1-eta}
	\eta(z)=\kappa/z ,
\end{equation}
which itself is equal to the inverse $-\kappa/z$ of the canonical link up to the minus sign.
In a sample $Y_1,\ldots,Y_n$ with covariates $\mathbf{x}_1,\ldots,\mathbf{x}_n$ the responses $Y_i$ are gamma distributed with means $\mu_i=\eta(\mathbf{f}(\mathbf{x}_i)^\trp\BS{\beta})$ and common shape parameter $\kappa$.

In an experimental design setup the covariates $\mathbf{x}_i$ may be chosen by the experimenter from an experimental region $\mathcal{X}$ over which the generalized linear model is assumed to be valid.
In particular, for gamma distributed responses the means $\mu_i$ have to be positive ($\mu(\mathbf{x}_i;\BS{\beta})>0$).
This implies the natural restriction on the parameter region $\mathcal{B}$ of potential values for the parameter vector $\BS{\beta}$ that for every $\BS{\beta}\in\mathcal{B}$ the linear component has to be positive ($\mathbf{f}(\mathbf{x})^\trp\BS{\beta}>0$) for all $\mathbf{x}\in\mathcal{X}$.
Further note that for reasons of parameter identifiability the regression functions $f_0,\ldots,f_{p-1}$ are assumed to be linearly independent on the experimental region $\mathcal{X}$. 

The aim of experimental design is to optimize the performance of the statistical analysis.
The contribution of an observation $Y_i$ to the performance is measured in terms of its information.
In the present generalized linear models framework, for a single observation at an experimental setting $\mathbf{x}$ the elemental information matrix is given by
\begin{equation}
\mathbf{M}(\mathbf{x};\BS{\beta})
=\lambda(\mathbf{f}(\mathbf{x})^\trp\BS{\beta})\,\mathbf{f}(\mathbf{x})\,\mathbf{f}(\mathbf{x})^\trp 
\label{eq2-4}
\end{equation} 
(see \citet{fedorov2013optimal} or  \citet{atkinson2015designs}) ,   
where $\lambda$ is a positive valued function which is called the intensity function. 
Note that through the intensity function the elemental information depends on the parameter vector $\BS{\beta}$.

In generalized linear models the intensity is given by
\begin{equation}
	\lambda(\mathbf{f}(\mathbf{x})^\trp\BS{\beta}) = \eta'(\mathbf{f}(\mathbf{x})^\trp\BS{\beta})^{2}/\mathrm{Var}(Y) .
	\label{lambda-glm}
\end{equation} 
In the case of a canonical link we have $\mathrm{Var}(Y) = \eta'(\mathbf{f}(\mathbf{x})^\trp\BS{\beta})$ and the intensity reduces to the variance.
In particular, in the gamma model with inverse link the intensity function is
\begin{equation}
\lambda(z)=\kappa/z^{2} ,
\label{lambda-gamma}
\end{equation} 
because the minus sign in the inverse of the link function does not affect the intensity
(cf.\ \citet{GAFFKE2019}).
The (per experiment) Fisher information of $n$ independent observations $Y_i$ at experimental settings $\mathbf{x}_i$ is then given by
 \begin{equation}
 \mathbf{M}(\mathbf{x}_1,\dots,\mathbf{x}_n;\BS{\beta})
 =\sum_{i=1}^n \mathbf{M}(\mathbf{x}_i;\BS{\beta})
 =\sum_{i=1}^n \lambda(\mathbf{f}(\mathbf{x}_i)^\trp\BS{\beta})\, \mathbf{f}(\mathbf{x}_i)\,\mathbf{f}(\mathbf{x}_i)^\trp .
 \label{eq2-4-Fisher}
\end{equation} 
The aim of finding an exact optimal design $\mathbf{x}_1^*,\ldots,\mathbf{x}_n^*$ is to optimize the Fisher information in a certain sense because their inverse is proportional to the asymptotic covariance matrix of the maximum likelihood estimator for $\BS{\beta}$ (see \citet{fahrmeir1985}).
 
 As this discrete optimization problem is too difficult, in general, we will deal with  approximate (continuous) designs $\xi$ in the spirit of \citet{kiefer1959jrssb} (see also \citet{silvey1980optimal}, p.15) throughout the remainder of the present paper.
 An approximate design $\xi$ is defined on the experimental region $\mathcal{X}$ by mutually distinct support points $\mathbf{x}_1,\ldots,\mathbf{x}_m$ and  corresponding weights $w_1,\ldots, w_m>0$ such that  $\sum_{i=1}^{m} w_i=1$. 
 In terms of an exact design the support points $\mathbf{x}_i$ may be interpreted as the distinct experimental settings and the weights $w_i$ as their corresponding relative frequencies in the sample.
 The relaxation of an approximate design is then that the weights $w_i$ may be chosen continuously and need not to be multiples of $1/n$.
The standardized (per observation) information matrix of a design $\xi$ is defined by 
\begin{equation}
	\mathbf{M}(\xi;\BS{\beta})
	=\sum_{i=1}^m w_i \mathbf{M}(\mathbf{x}_i;\BS{\beta})
	=\sum_{i=1}^m w_i \lambda(\mathbf{f}(\mathbf{x}_i)^\trp\BS{\beta})\, \mathbf{f}(\mathbf{x}_i)\,\mathbf{f}(\mathbf{x}_i)^\trp .
	\label{eq2-4-approximate}
\end{equation} 

Design optimization is now concerned with finding an approximate design $\xi^*$ which maximizes the Fisher information or, equivalently, minimizes the asymptotic covariance matrix. 
Because matrices cannot be compared directly, in general, some convex real-valued criterion function $\Phi$ to be minimized will be considered instead which commonly depends on the design $\xi$ through the Fisher information $\mathbf{M}(\xi; \BS{\beta})$.
A design $\xi^*$  will then be called $\Phi$-optimal when it minimizes $\Phi(\xi)$, $\Phi(\xi^*)=\min \Phi(\xi)$.
As the information matrix depends on the parameter vector $\BS{\beta}$ the so obtained designs $\xi^*$ are locally $\Phi$-optimal at a given parameter value $\BS{\beta}$ (\citet{chernoff1953}) and may change with $\BS{\beta}$. 
To avoid the parameter dependence so-called robust versions of the criteria can be considered like ``Bayesian'' criteria which involve a weighting measure (``prior'') on the parameters (see \citet{AtkinsonDonevTobias}, ch.~18) or ``minimax'' criteria which aim at minimizing the worst case scenario for the parameter settings (see the ``standardized minimax'' criteria in \citet{10.2307/2345917}).
In the following we will focus on the local \textit{D}- and \textit{IMSE}-criteria and the corresponding maximin efficiency (``standardized maximin'') criteria. 

The \textit{D}-criterion is the most commonly used design criterion.
It is related to the estimation of the model parameters $\BS{\beta}$ and aims at minimizing the determinant of the asymptotic covariance matrix, $\Phi(\mathbf{M})=\det(\mathbf{M}^{-1})$ for positive definite information matrix $\mathbf{M}$, and $\Phi(\mathbf{M})=\infty$ for singular $\mathbf{M}$. 
A design $\xi^*$ is then called locally \textit{D}-optimal at $\BS{\beta}$ when $\det(\mathbf{M}(\xi^*; \BS{\beta})^{-1})=\min\det(\mathbf{M}(\xi; \BS{\beta})^{-1})$.
The \textit{D}-criterion can be motivated by the fact that it measures the (squared) volume of the asymptotic confidence ellipsoid of the maximum likelihood estimator for $\BS{\beta}$.
However, its popularity predominantly stems from its nice analytic properties.

Note that in the present situation the property of $\mathbf{M}(\xi; \BS{\beta})$ being nonsingular does not depend on the value of the parameter vector $\BS{\beta}$ because the intensity $\lambda(\mathbf{f}(\mathbf{x})^\trp\BS{\beta})$ is greater than zero for all $\mathbf{x}\in \mathcal{X}$ and all $\BS{\beta}\in \mathcal{B}$. 

The definition of the \textit{IMSE}-criterion is based on the estimation (prediction) of the mean response $\mu(\mathbf{x};\BS{\beta})$.
It aims at minimizing the average asymptotic variance of the predicted mean response $\hat\mu(\mathbf{x})=\mu(\mathbf{x};\hat{\BS{\beta}})$, where averaging is taken with respect to a standardized measure $\nu$ on $\mathcal{X}$ (see \citet{li2018optimal} and \citet{doi:10.1002/cjs.11571}). 
For a generalized linear model the asymptotic variance is given by 
\begin{equation}
	\mathrm{asVar}(\hat\mu(\mathbf{x});\xi,\BS{\beta})
	=\eta'(\mathbf{f}(\mathbf{x})^\trp\BS{\beta})^{2}
	  \mathbf{f}(\mathbf{x})^\trp\mathbf{M}(\xi; \BS{\beta})^{-1} \mathbf{f}(\mathbf{x}),
\label{as-var-mean-response}
\end{equation}
for all $\mathbf{x}\in \mathcal{X}$. 
For a canonical link we have $\lambda=\eta'$ and hence
\begin{equation}
	\mathrm{asVar}(\hat\mu(\mathbf{x});\xi,\BS{\beta})
	=\lambda(\mathbf{f}(\mathbf{x})^\trp\BS{\beta})^{2}
	\mathbf{f}(\mathbf{x})^\trp\mathbf{M}(\xi; \BS{\beta})^{-1} \mathbf{f}(\mathbf{x}).
	\label{as-var-mean-response-canonical-link}
	\end{equation}
The integrated mean-squared error (IMSE) is then defined as the average prediction variance
\begin{eqnarray}
\mathrm{IMSE}(\xi;\BS{\beta},\nu)
&=&
\int \mathrm{asVar}(\hat\mu(\mathbf{x});\xi,\BS{\beta})\,\nu(\mathrm{d} \mathbf{x})
\nonumber
\\
&=&
\int \lambda(\mathbf{f}(\mathbf{x})^\trp\BS{\beta})^{2}
\mathbf{f}(\mathbf{x})^\trp\mathbf{M}(\xi; \BS{\beta})^{-1} \mathbf{f}(\mathbf{x})\,
\nu(\mathrm{d} \mathbf{x})
\label{IMSE}
\end{eqnarray}
with respect to a given standardized measure $\nu$ on the experimental region $\mathcal{X}$ ($\nu(\mathcal{X})=1$).

By a standard method to express the \textit{IMSE}-criterion (see e.\,g.\  \citet{doi:10.1002/cjs.11571}) the asymptotic variance can be rewritten as
\begin{equation}
	\mathrm{asVar}(\hat\mu(\mathbf{x});\xi,\BS{\beta})
	=\mathrm{trace}(\lambda(\mathbf{f}(\mathbf{x})^\trp\BS{\beta})^{2}
	\mathbf{f}(\mathbf{x})\mathbf{f}(\mathbf{x})^\trp\mathbf{M}(\xi; \BS{\beta})^{-1}) .
	\label{as-var-mean-response-trace}
\end{equation}
Hence, the IMSE is given by
\begin{equation}
\mathrm{IMSE}(\xi;\BS{\beta},\nu)
=
\mathrm{trace}(\mathbf{V}(\BS{\beta};\nu)\,\mathbf{M}(\xi; \BS{\beta})^{-1}) ,
\label{IMSE-trace}
\end{equation}
where 
\begin{equation}
\mathbf{V}(\BS{\beta};\nu)
=\int \lambda(\mathbf{f}(\mathbf{x})^\trp\BS{\beta})^{2} 
  \mathbf{f}(\mathbf{x})\mathbf{f}(\mathbf{x})^\trp\,\,\nu(\mathrm{d} \mathbf{x})  
  \label{mom-mat}
\end{equation}
denotes a weighted ``moment'' matrix with respect to the measure $\nu$. 
Note that the leading term under the integral in $\mathbf{V}(\BS{\beta};\nu)$ differs from that in the virtual information matrix $\mathbf{M}(\nu; \BS{\beta})$ by replacing the intensity $\lambda$ by $\lambda^2$.
Moreover, in contrast to the \textit{D}-criterion, the \textit{IMSE}-criterion does not solely depend on the information matrix $\mathbf{M}(\xi;\BS{\beta})$, but also depends through the weighting matrix $\mathbf{V}(\BS{\beta};\nu)$ explicitly on the parameter vector $\BS{\beta}$ and additionally on the measure $\nu$ as a supplementary argument.
The \textit{IMSE}-criterion is thus defined by
$\Phi(\mathbf{M};\BS{\beta},\nu)=\mathrm{trace}(\mathbf{V}(\BS{\beta};\nu)\,\mathbf{M}^{-1})$.
A design $\xi^*$ is then called locally \textit{IMSE}-optimal with respect to $\nu$ at $\BS{\beta}$ when $\mathrm{trace}(\mathbf{V}(\BS{\beta};\nu)\mathbf{M}(\xi^*; \BS{\beta})^{-1})=\min\mathrm{trace}(\mathbf{V}(\BS{\beta};\nu)\mathbf{M}(\xi; \BS{\beta})^{-1})$.

To avoid the parameter dependence of an optimal design under local criteria, we will also consider as a ``robust'' alternatives maximin efficiency criteria which are also called standardized optimality criteria (see \citet{DETTE20064397}).
For this we first have to introduce the concept of efficiency.
Let the local criterion $\Phi_{\BS{\beta}}$ at $\BS{\beta}$ depend homogeneously on the information matrix, i.\,e.,\ $\Phi_{\BS{\beta}}(\xi)=\phi(\mathbf{M}(\xi;\BS{\beta}))$ for some function $\phi$ on the set of positive definite matrices satisfying $\phi(c\mathbf{M})=c^{-1}\phi(\mathbf{M})$ for $c>0$ (cf \citet{pukelsheim1993optimal}, ch.~5, for the related concept of information functions).
Then the efficiency of a design $\xi$ (locally at $\BS{\beta}$) is defined by
\begin{equation*}
	\mathrm{eff}(\xi;\BS{\beta})=\frac{\Phi_{\BS{\beta}}(\xi_{\BS{\beta}}^*)}{\Phi_{\BS{\beta}}(\xi)},
	\label{def-eff}
\end{equation*}
where $\xi_{\BS{\beta}}^*$ is the $\Phi$-optimal design (locally at $\BS{\beta}$).
The efficiency can then be interpreted as the proportion of observations required under the optimal design $\xi_{\BS{\beta}}^*$ to obtain the same value of $\Phi$ as for design $\xi$.
For example, an efficiency of $0.5$ means that with an optimum design $\xi_{\BS{\beta}}^*$ only half as many observations as for $\xi$ are necessary to get the same precision. 
Maximin efficiency then aims at maximizing the worst efficiency $\inf_{\BS{\beta}\in\mathcal{B}^{\prime}} \mathrm{eff}_{\Phi}(\xi;\BS{\beta})$ over a given subset $\mathcal{B}^{\prime}$ of interest of the parameter region $\mathcal{B}$.
In order to arrive at a minimization problem we define the maximin efficiency criterion by the inverse relation
\begin{equation}
	\Phi(\xi)=\sup_{\BS{\beta}\in\mathcal{B}^{\prime}}\frac{\Phi_{\BS{\beta}}(\xi)}{\Phi_{\BS{\beta}}(\xi_{\BS{\beta}}^*)} .
	\label{def-maximin-eff}
\end{equation}
Note that $\Phi$ is convex if the local criteria $\Phi_{\BS{\beta}}$ are all convex. 

For maximin \textit{D}-efficiency we have to choose the homogeneous version
$\Phi_{\BS{\beta}}(\xi)=(\det(\mathbf{M}(\xi;\BS{\beta})))^{-1/p}$ of the local \textit{D}-criterion (see \citet{pukelsheim1993optimal}, ch.~6) to get the maximin \textit{D}-efficiency criterion 
\[
	\Phi_{\textit{D-ME}}(\xi)=
	\sup_{\BS{\beta}\in\mathcal{B}^{\prime}}
	\left(
	\frac{\det(\mathbf{M}(\xi;\BS{\beta}))}{\det(\mathbf{M}(\xi_{\BS{\beta}}^*;\BS{\beta}))}
	\right)^{-1/p} ,
\]
where $\xi_{\BS{\beta}}^*$ denotes the locally \textit{D}-optimal design at $\BS{\beta}$.
A design $\xi^*$ is then called maximin \textit{D}-efficient on $\mathcal{B}^{\prime}$ when $\Phi_{\textit{D-ME}}(\xi^*)=\min\Phi_{\textit{D-ME}}(\xi)$.

The local \textit{IMSE}-criterion is already homogeneous because it is a linear criterion.
Thus the maximin \textit{IMSE}-efficiency criterion can be defined directly as
\[
\Phi_{\textit{IMSE-ME}}(\xi;\nu)=\sup_{\BS{\beta}\in\mathcal{B}^{\prime}}\frac{\mathrm{IMSE}(\xi,\BS{\beta},\nu)}{\mathrm{IMSE}(\xi^*_{\BS{\beta}},\BS{\beta},\nu)},
\]
where $\xi_{\BS{\beta}}^*$ denotes the locally \textit{IMSE}-optimal design at $\BS{\beta}$.
A design $\xi^*$ is then called maximin \textit{IMSE}-efficient with respect to $\nu$ on $\mathcal{B}^{\prime}$ when $\Phi_{\textit{IMSE-ME}}(\xi^*;\nu)=\min\Phi_{\textit{IMSE-ME}}(\xi;\nu)$.

In particular, for the gamma model with inverse link we have $\lambda(z)=\kappa/z^2$ (see (\ref{lambda-gamma})) which implies that
\begin{equation}
	\mathbf{M}(\xi;\BS{\beta})
	=\sum_{i=1}^m w_i \kappa(\mathbf{f}(\mathbf{x}_i)^\trp\BS{\beta})^{-2}\, \mathbf{f}(\mathbf{x}_i)\,\mathbf{f}(\mathbf{x}_i)^\trp 
	\label{eq2-4-approximate-gamma}
\end{equation} 
and
\begin{equation}
	\mathbf{V}(\BS{\beta};\nu)
	=\int \kappa^2(\mathbf{f}(\mathbf{x})^\trp\BS{\beta})^{-4} 
	\mathbf{f}(\mathbf{x})\mathbf{f}(\mathbf{x})^\trp\,\,\nu(\mathrm{d} \mathbf{x})  .
	\label{mom-mat-gamma}
\end{equation}
Hence, in both the \textit{D}- and the \textit{IMSE}-criterion the shape parameter $\kappa$ occurs only as a factor which does not affect the optimization problem. 
Without loss of generality we may thus assume $\kappa=1$ in the remainder of the text.

\section{Equivariance}
\label{sec-3}

In- and equivariance play an important role for optimal design in linear models.
However, these concepts can also be applied in the context of generalized linear models as established in \citet{radloff2016invariance}.

The essential idea of equivariance in the design setup is to transfer an already known optimal design on a given (standardized) experimental region to another experimental region of interest by a suitable transformation while keeping the model structure unchanged.
The most prominent approach of this kind is the method of canonical transformation propagated by \citet{10.2307/2346142}.

Throughout we accompany each conceptual step by a simple running example (Example~\ref{ex-1}).
We start with a one-to-one transformation   
$g:\,\mathcal{X}\to \mathcal{Z}$
which maps the experimental region $\mathcal{X}$ onto a potentially different region $\mathcal{Z}$.

\begin{exmp}
\label{ex-1}
Let $\mathcal{X}=[0,1]$ be the one-dimensional standard unit interval and $\mathcal{Z}=[a,b]$ another non-degenerate interval, $b>a$. 
Then the shift and scale transformation $g(x)=a+cx$, where $c=b-a$, maps $\mathcal{X}$ onto  $\mathcal{Z}$.
\eop
\end{exmp}

The next ingredient connects the transformation $g$ with the vector of regression functions:
$\mathbf{f}$ is said to be linearly equivariant with respect to $g$ if there exists 
a (nonsingular) matrix  $\mathbf{Q}_g$ such that 
$\mathbf{f}(g(\mathbf{x}))=\mathbf{Q}_g\mathbf{f}(\mathbf{x})$ for all $\mathbf{x}\in \mathcal{X}$, which will be assumed to hold throughout the remainder of this text.

\begin{exmp*}[Example~\ref{ex-1} continued]
Let $\mathbf{f}(x)=(1,x)^\trp$ the vector of regression functions for a simple one-dimensional linear regression, $p=2$, such that the linear component is $\mathbf{f}(x)^\trp\BS{\beta}=\beta_0+\beta_1 x$. Then for $g(x)=a+cx$ the transformation matrix $\mathbf{Q}_g$ is given by 
\[
\mathbf{Q}_g =
\left(
\begin{array}{cc}
	1 & 0
	\\
	a & c
\end{array}
\right).
\]
\eop
\end{exmp*}

In contrast to the situation in linear models, additionally a transformation $\tilde{g}:\,\mathcal{B}\to\tilde{\mathcal{B}}$ of the parameter vector $\BS{\beta}$ is required in the present setup of generalized linear models.
This approach of equivariance with respect to a pair $(g,\tilde{g})$ of transformations of the settings $\mathbf{x}$ and the parameters $\BS{\beta}$, respectively, is in accordance with the general concept of equivariance in statistical analysis (see e.\,g.\ \citet{lehmann1959testing}, ch.~6).

A natural choice for the transformation $\tilde{g}$ is a reparameterization which leaves the value of the linear component unchanged, $\mathbf{f}(g(\mathbf{x}))^\trp\tilde{g}(\BS{\beta})=\mathbf{f}(\mathbf{x})^\trp\BS{\beta}$ for all $\mathbf{x}\in\mathcal{X}$.
This is accomplished by setting $\tilde{g}(\BS{\beta})=\mathbf{Q}_g^{-\trp}\BS{\beta}$, where ``$\cdot^{-\trp}$'' denotes the inverse of a transposed matrix, and $\tilde{\mathcal{B}}=\tilde{g}(\mathcal{B})$.
For convenience we denote $\tilde{g}(\BS{\beta})$ by $\tilde{\BS{\beta}}=(\tilde{\beta}_0,\ldots, \tilde{\beta}_{p-1})^\trp$.  

\begin{exmp*}[Example~\ref{ex-1} continued]
	For  $g(x)=a+cx$ and simple linear regression $\mathbf{f}(x)=(1,x)^\trp$ the transformation matrix for the parameter vector is 
	\[
	\mathbf{Q}_g^{-\trp} =
	\left(
	\begin{array}{cc}
		1 & -a/c
		\\
		0 & 1/c
	\end{array}
	\right),
	\]
	and the transformation $\tilde{g}(\BS{\beta})=\mathbf{Q}_g^{-\trp}\BS{\beta}$ results in $\tilde{\beta_0}=\beta_0-a\beta_1/c$ and $\tilde{\beta_1}=\beta_1/c$.
	If for given value of $\BS{\beta}=(\beta_0,\beta_1)^\trp$ the pair $(g,\tilde{g})$ is chosen in such a way that $\tilde{\BS{\beta}}=(0,1)^\trp$, i.\,e.,\ $c=\beta_1$ and $a=\beta_0$, then $g$ represents essentially the canonical transformation used in \citet{10.2307/2346142}.
	
	For the gamma model the parameter region $\mathcal{B}$ is restricted by the constraint that the linear component $\mathbf{f}(x)^\trp\BS{\beta}=\beta_0+\beta_1 x$ is positive for all $x\in\mathcal{X}=[0,1]$.
	Hence, the maximal parameter region is $\mathcal{B}=\{\BS{\beta};\,\beta_0>0, \beta_1>-\beta_0\}$ which is displayed in Figure~\ref{fig:opt1}.
	The transformed parameter region is then $\tilde{\mathcal{B}}=\tilde{g}(\mathcal{B})=\{\tilde{\BS{\beta}};\,\tilde{\beta_0}+\tilde{\beta_1}a>0, \tilde{\beta_0}+\tilde{\beta_1}b>0\}$.
	In particular, to obtain the symmetric unit interval $\mathcal{Z}=[-1,1]$ as secondary experimental region, the transformation $g(x)=2x-1$ is to be chosen with $a=-1$ and $c=2$, and the transformed parameter region becomes $\tilde{\mathcal{B}}=\{\tilde{\BS{\beta}};\,|\tilde{\beta_1}|<\tilde{\beta_0}\}$.
    \begin{figure}[H]
	\centering
	\includegraphics[width=9cm, height=6.5cm]{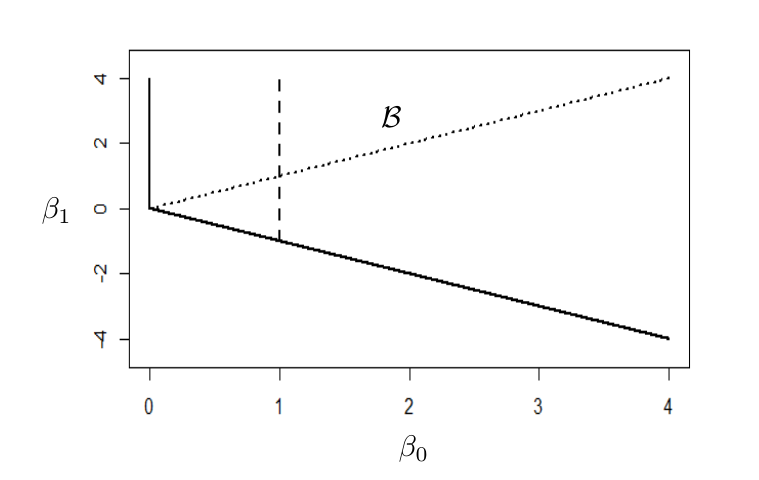}
	\caption{Parameter region $\mathcal{B}$ for the one-factor gamma model on $[0,1]$}
	\label{fig:opt1}
\end{figure}
\eop
\end{exmp*}

Note that for each pair $(g, \tilde{g})$ of transformations the mean response and the intensity remain unchanged,  $\mu(g(\mathbf{x});\tilde{g}(\BS{\beta}))=\mu(\mathbf{x};\BS{\beta})$ and $\lambda(\mathbf{f}(g(\mathbf{x}))^\trp\tilde{g}(\BS{\beta}))=\lambda(\mathbf{f}(\mathbf{x})^\trp\BS{\beta})$.
Having this in mind we study how these transformations act on a design and its information matrix:
For a design $\xi$ with support points $\mathbf{x}_i$ and corresponding weights $w_i$, $i=1,\ldots,m$, we denote by $\xi^g$ its image under the transformation $g$, i.\,e.,\ $\xi^g$ has support points $\mathbf{z}_i=g(\mathbf{x}_i)$ with weights $w_i$, $i=1,\ldots,m$, respectively, and is hence a design on $\mathcal{Z}$.
Then for the associated information matrices we obtain
\begin{eqnarray}
	\mathbf{M}(\xi^{g};\tilde{g}(\BS{\beta}))
	&=&
	  \sum_{i=1}^m w_i  \lambda(\mathbf{f}(g(\mathbf{x}_i))^\trp\tilde{g}(\BS{\beta})) \mathbf{f}(g(\mathbf{x}_i))\mathbf{f}(g(\mathbf{x}_i))^\trp
	 \nonumber
	 \\
	&=&
	 \sum_{i=1}^m w_i  \lambda(\mathbf{f}(\mathbf{x}_i)^\trp\BS{\beta}) \mathbf{Q}_{g}\mathbf{f}(\mathbf{x}_i)\mathbf{f}(\mathbf{x}_i)^\trp\mathbf{Q}^\trp_{g}
	 \nonumber
	\\
	&=&
	  \mathbf{Q}_{g} \left(\sum_{i=1}^m w_i \lambda(\mathbf{f}(\mathbf{x}_i)^\trp\BS{\beta}) \mathbf{f}(\mathbf{x}_i)\mathbf{f}(\mathbf{x}_i)^\trp\right) \mathbf{Q}^\trp_{g}
	 = \mathbf{Q}_{g}\mathbf{M}(\xi; \BS{\beta})\mathbf{Q}^\trp_{g} 
	 \label{info-transformation}
\end{eqnarray} 
(see \citet{radloff2016invariance}).
In short, the pair of simultaneous transformations  $(g, \tilde{g})$ induces the transformation $\mathbf{M}(\xi;\BS{\beta})\to\mathbf{Q}_g\mathbf{M}(\xi;\BS{\beta})\mathbf{Q}_g^\trp$ of the  information matrix.

\begin{exmp*}[Example~\ref{ex-1} continued]
	Let $\xi$ be supported on the endpoints $x_1=0$ and $x_2=1$ of the experimental region $\mathcal{X}=[0,1]$ with corresponding weights $w_1=1-w$ and $w_2=w$, respectively. 
	 For the gamma model with simple linear regression, $\mathbf{f}(x)=(1,x)^\trp$, denote by $\lambda_0=\lambda(\beta_0)$ and $\lambda_1=\lambda(\beta_0+\beta_1)$ the intensities at the support points $0$ and $1$.
	The information matrix of $\xi$ is given by
	\[
		\mathbf{M}(\xi;\BS{\beta})
		=
		\left(
		\begin{array}{cc}
			(1-w)\lambda_0 + w\lambda_1 & w\lambda_1
			\\
			w\lambda_1 & w\lambda_1
		\end{array}
		\right).
	\]
	For $g(x)=a+cx$ the induced design $\xi^g$ is supported on the endpoints $z_1=a$ and $z_2=b$ of the induced experimental region $\mathcal{Z}=[a,b]$ with weights $1-w$ at $a$ and $w$ at $b$. 
	Under $\tilde{\BS{\beta}}=\tilde{g}(\BS{\beta})$ the intensities at $a$ and $b$ are $\lambda_0$ and $\lambda_1$, respectively, and the information matrix of $\xi^g$ is
	\[
		\mathbf{M}(\xi^g;\tilde{g}(\BS{\beta}))
		=
		\left(\begin{array}{cc}
			(1-w)\lambda_0 + w\lambda_1 & (1-w)\lambda_0 a + w\lambda_1 b
			\\
			(1-w)\lambda_0 a + w\lambda_1 b & (1-w)\lambda_0 a^2 + w\lambda_1 b^2
		\end{array}\right)
		=
		\mathbf{Q}_g \mathbf{M}(\xi;\BS{\beta}) \mathbf{Q}_g^\trp.
	\]
\eop
\end{exmp*}

The final step is the equivariance of the criterion $\Phi$.
In analogy to the terminology in \citet{HEILIGERS1992113} for linear models we will call a convex optimality criterion $\Phi$ equivariant with respect to a transformation $g$, if $\Phi$ preserves the ordering under the transformation $g$, i.\,e.,\ for any two designs $\xi_1$ and $\xi_2$ the relation $\Phi(\xi_1)\leq \Phi(\xi_2)$ implies $\Phi(\xi_1^g)\leq \Phi(\xi_2^g)$. 

As in the present situation of generalized linear models also the parameter vector $\BS{\beta}$ and eventually some supplementary arguments have to be changed in the criterion during the transformation, more care has to be taken.
We therefore introduce a second criterion function $\Phi^{\prime}=\Phi_{g,\tilde{g}}$ for the designs on $\mathcal{Z}$ which may depend on the transformations $g$ and $\tilde{g}$.
Then we will call a pair of criteria $\Phi$ and $\Phi^{\prime}$ equivariant with respect to the pair  $(g,\tilde{g})$ of transformations, when the ordering is preserved, i.\,e.,\ the relation $\Phi(\xi_1)\leq  \Phi(\xi_2)$ implies $\Phi^{\prime}(\xi_1^g)\leq \Phi^{\prime}(\xi_2^g)$. 

With this definitions we obtain the following result that in the case of equivariance the optimality of designs is preserved under transformations.
\begin{theorem}
	\label{th-equivariance}
	Let the pair of criteria $\Phi$ and $\Phi^{\prime}$ be equivariant with respect to the pair  $(g,\tilde{g})$ of transformations.
	If $\xi^*$ is $\Phi$-optimal, then its image $(\xi^*)^g$ is $\Phi^{\prime}$-optimal.
\end{theorem}

We will now establish that the \textit{D}- and \textit{IMSE}-criterion are equivariant, if simultaneously the parameter vector $\BS{\beta}$ and eventual supplementary arguments are transformed.
By (\ref{info-transformation}) we obtain for the \textit{D}-criterion
\begin{equation}
	\label{eq-transformation-d}
	\det(\mathbf{M}(\xi^g;\tilde{g}(\BS{\beta}))^{-1}) = \det(\mathbf{Q}_g)^{-2}\det(\mathbf{M}(\xi;\BS{\beta})^{-1}).
\end{equation}
Let $\Phi$ be the local \textit{D}-criterion at $\BS{\beta}$ and $\Phi^{\prime}$ the local \textit{D}-criterion at $\tilde{g}(\BS{\beta)}$, then the \textit{D}-criterion is equivariant under simultaneous transformation of $\BS{\beta}$ and by Theorem~\ref{th-equivariance} the locally \textit{D}-optimal design can be transferred.
\begin{corollary}
	\label{cor-equivariance-d}
	If $\xi^*$ is locally \textit{D}-optimal on $\mathcal{X}$ at $\BS{\beta}$, then $(\xi^*)^g$  is locally \textit{D}-optimal on $\mathcal{Z}$ at $\tilde{\BS{\beta}}=\tilde{g}(\BS{\beta})$.
\end{corollary}

\begin{exmp*}[Example~\ref{ex-1} continued]
	For the gamma model with simple linear regression, $\mathbf{f}(x)=(1,x)^\trp$, the locally \textit{D}-optimal design $\xi^{*}$ on the unit interval $\mathcal{X}=[0,1]$ is supported by the endpoints $x_1=0$ and $x_2=1$ and assigns equal weights $w^*=1/2$ to these endpoints for any value of the parameter vector $\BS{\beta}\in\mathcal{B}$ (see \citet{GAFFKE2019}).
	Then for any other interval $\mathcal{Z}=[a,b]$ as the experimental region we may consider the transformation $g(x)=a+cx$, $c=b-a$, together with $\tilde{g}(\BS{\beta})=\mathbf{Q}_g^{-\trp}\BS{\beta}$. 
	By Corollary~\ref{cor-equivariance-d} the design $(\xi^*)^g$ which assigns equal weights $w^*=1/2$ to he endpoints $z_1=a$ and $z_2=b$ of the experimental region $\mathcal{Z}$ is locally \textit{D}-optimal for any value of the parameter vector $\tilde{\BS{\beta}}=\tilde{g}(\BS{\beta})\in\tilde{\mathcal{B}}=\tilde{g}(\mathcal{B})$.
\eop
\end{exmp*}

In the situation of Example~\ref{ex-1} the locally \textit{D}-optimal design does not depend on the parameter $\BS{\beta}$. This will typically not hold true, if the underlying model for the linear component becomes more complex.

\begin{figure}[h]
 \centering
 \scalebox{0.9}{
\begin{minipage}[b]{0.5\textwidth}
    \includegraphics[width=7cm, height=6cm]{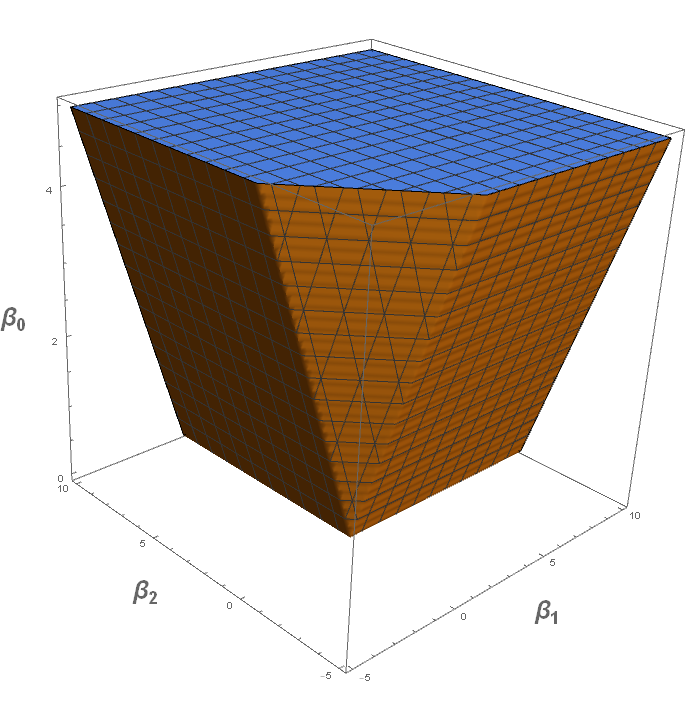}
  \end{minipage}}
\hfill
\scalebox{0.9}{
  \begin{minipage}[b]{0.5\textwidth}
    \includegraphics[width=9cm, height=6.5cm]{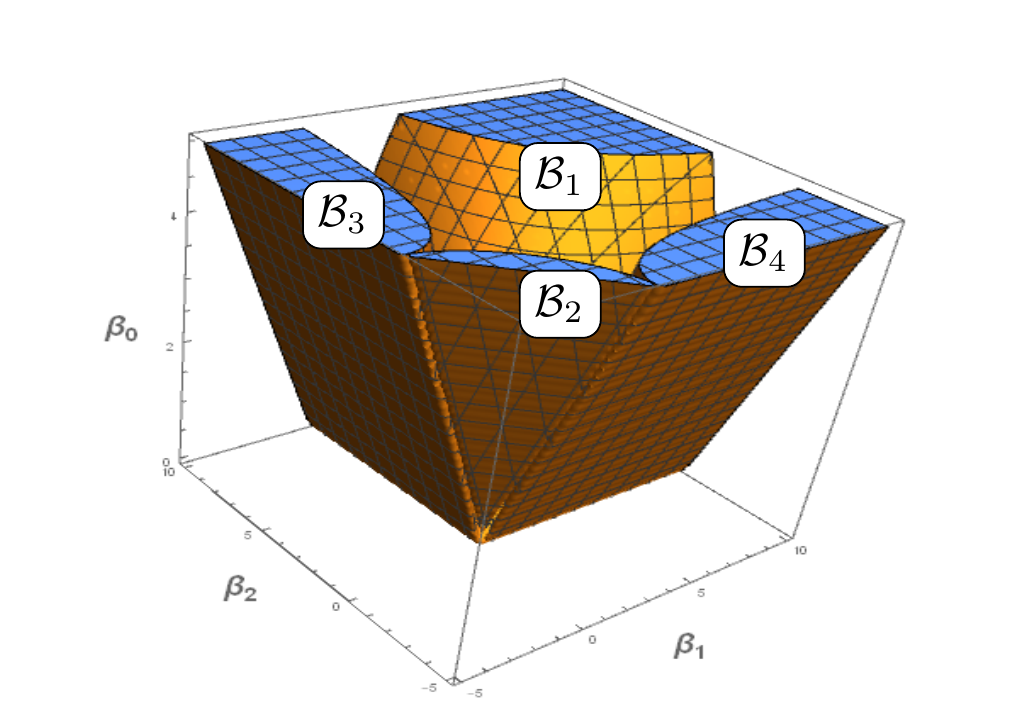}
  \end{minipage}
  } \caption{Parameter region $\mathcal{B}$ in the two-factor gamma model on $[0,1]^2$ (left panel); subregions $\mathcal{B}_k$, $k=1,\ldots,4$ of parameters with minimally supported locally \textit{D}-optimal designs (right panel)}
 \label{fig:opt2}
\end{figure}

\begin{exmp}
	\label{exmp1}
 We consider the gamma model with the linear component $\mathbf{f}(\mathbf{x})^\trp\BS{\beta}=\beta_0+\beta_1 x_1+\beta_2 x_2$ of a multiple linear regression of two covariates, $\mathbf{x}=(x_1,x_2)^\trp$, where $\mathbf{f}(\mathbf{x})=(1,x_1,x_2)^\trp$, $p=3$, with  the unit square $\mathcal{X}=[0,1]^2$
as experimental region.
Denote by $\mathbf{x}_1=(0,0)^\trp$,  $\mathbf{x}_2=(1,0)^\trp$, $\mathbf{x}_3=(0,1)^\trp$ and $\mathbf{x}_4=(1,1)^\trp$ the vertices of $\mathcal{X}$. 
The parameter region $\mathcal{B}$ is the set of all parameter vectors $\BS{\beta}=(\beta_0,\beta_1,\beta_2)^\trp$  such that the linear component at $\mathbf{x}_1,\ldots,\mathbf{x}_4$ is positive, i.\,e.,\ $\beta_0>0$, $\beta_0+\beta_1>0$, $\beta_0+\beta_2>0$, and $\beta_0+\beta_1+\beta_2>0$.
This region, depicted in the left panel of Figure \ref{fig:opt2}, constitutes a cone in the three-dimensional Euclidean space.

According to \citet{10.2307/2336960} the minimally supported design $\xi^*$ which assigns equal weights $w_i^*=1/3$ to the support points $\mathbf{x}_i$,  $i=1,2,3$, is locally \textit{D}-optimal at $\BS{\beta}$, when $\BS{\beta}$ satisfies $\beta_0^2-\beta_1\beta_2\leq 0$.  
The subset $\mathcal{B}_1$ of these $\BS{\beta}$ in $\mathcal{B}_1$ is shown in the right panel of Figure \ref{fig:opt2}.  

Now equivariance can be used to find \textit{D}-optimal designs for other parameter values different from those in $\mathcal{B}_1$.
For this we use transformations which map the experimental region onto itself, $\mathcal{Z}=\mathcal{X}$:
\begin{equation}
	g_2(\mathbf{x})=(1-x_1,1-x_2)^\trp,
	\quad
	g_3(\mathbf{x})=(1-x_1,x_2)^\trp 
	\quad \mathrm{ and } \quad
	g_4(\mathbf{x})=(x_1,1-x_2)^\trp.
	\label{exp-trans}
\end{equation}
Here $g_3$ and $g_4$ represent the reflection with respect to the first and second covariate $x_1$ and $x_2$, respectively, and $g_2$ is the simultaneous reflection with respect to both covariates. 
Alternatively, $g_2$ can also be described as a rotation by $180$ degree.
We also introduce $g_1=\mathrm{id}$ as the identity mapping.

The regression function $\mathbf{f}(\mathbf{x})=(1,x_1,x_2)^\trp$ is linearly equivariant with respect to these transformations with corresponding matrices
\[
\mathbf{Q}_{g_2}=\left(\begin{array}{rrr}1&0&0\\1&-1&0\\1&0&-1\end{array}\right),
\,\,~~~\mathbf{Q}_{g_3}=\left(\begin{array}{rrr}1&0&0\\1&-1&0\\0&0&1\end{array}\right),\,\,
~~~\mathbf{Q}_{g_4}=\left(\begin{array}{rrr}1&0&0\\0&1&0\\1&0&-1\end{array}\right).
\]
For each $g_k$, $k=2,3,4$, the corresponding parameter transformation is given by $\tilde{g}_k(\BS{\beta})=\mathbf{Q}_{g_k}^{\trp}\BS{\beta}$.
Because $g_k$ maps the experimental region $\mathcal{X}$ onto itself, also the related transformation $\tilde{g}_k$ maps the parameter regions onto itself, $\tilde{\mathcal{B}}=\mathcal{B}$.

Starting from the parameter subregion $\mathcal{B}_1$, where the design $\xi^*$ is locally \textit{D}-optimal, we can define parameter subregions $\mathcal{B}_k=\tilde{g}_k(\mathcal{B}_1)$ induced by the transformations $g_k$, $k=2,3,4$. 
Explicit parameter formulations are given in Table~\ref{T-1} for these subregions, and they are also shown in the right panel of  Figure~\ref{fig:opt2}. 
All these subregions constitute cones.
Now by equivariance we can conclude that the designs $\xi^*_k=(\xi^*)^{g_k}$ are locally \textit{D}-optimal at $\BS{\beta}$ for $\BS{\beta}\in\mathcal{B}_k$.
The results are explicitly stated in Table~\ref{T-1}.

Note that the same optimal designs have been obtained before in \citet{idais2020local} by a straightforward application of the celebrated Kiefer-Wolfowitz equivalence theorem (see e.\,g.\ \citet{silvey1980optimal}).
Further note that the interior region shown in the right panel of Figure~ \ref{fig:opt2} contains those values for the parameter vector $\BS{\beta}$ for which a locally \textit{D}-optimal designs are supported on all four vertices and the corresponding weights depend on the values of $\BS{\beta}$ (see \citet{idais2020local}). 
\eop
\end{exmp}

\begin{table}[H]
\centering
{
\caption{Minimally supported locally \textit{D}-optimal designs and optimality regions for the two-factor gamma model on $[0,1]^2$. } 
\label{T-1}
\scalebox{0.9}{
\resizebox{\textwidth}{!}{
\begin{tabular}{lllll }
 \hline
\underline{ Transformation }&\,\,\, \underline{ Induced design }& \,\,\,\underline{ Induced support }&\,\,\,\,\,\underline{ Optimality region }&\,\,\,\,\,\,\,\,\,\underline{ Parameter formulation }  \\
\hspace{3ex}$g_1$             &\hspace{5ex}$\xi_1^{*}=\xi^*$ &\hspace{3ex}$\mathbf{x}_1$,  $\mathbf{x}_2$, $\mathbf{x}_3$&\hspace{3ex} $\mathcal{B}_1$ &$\beta_0^2- \beta_1\beta_2\leq 0$\\[0.5ex]
\hspace{3ex}$g_2$             &\hspace{5ex}$\xi_2^{*}=(\xi^*)^{g_2}$ &\hspace{3ex}$\mathbf{x}_2$,  $\mathbf{x}_3$, $\mathbf{x}_4$&\hspace{3ex} $\mathcal{B}_2=\tilde{g}_2(\mathcal{B}_1)$ &$(\beta_0+\beta_1+\beta_2)^2- \beta_1\beta_2\leq 0$\\[0.5ex]
\hspace{3ex}$g_3$           &\hspace{5ex}$\xi_3^{*}=(\xi^*)^{g_3}$ &\hspace{3ex}$\mathbf{x}_1$,  $\mathbf{x}_2$, $\mathbf{x}_4$&  \hspace{3ex}     $\mathcal{B}_3=\tilde{g}_3(\mathcal{B}_1)$ & $(\beta_0+\beta_1)^2+\beta_1\beta_2\leq 0$ \\[0.5ex]
\hspace{3ex}$g_4$           &\hspace{5ex}$\xi_4^{*}=(\xi^*)^{g_4}$   &\hspace{3ex}$\mathbf{x}_1$,  $\mathbf{x}_3$, $\mathbf{x}_4$&\hspace{3ex} 	$\mathcal{B}_4=\tilde{g}_4(\mathcal{B}_1)$ & $(\beta_0+\beta_2)^2+\beta_1\beta_2\leq 0$ \\
 \hline
\end{tabular}}}
}
\end{table}

Next we investigate equivariance for the \textit{IMSE}-criterion. 
There also the supplementary argument of the weighting measure $\nu$ has to be transformed.
Similar to the information matrix in (\ref{info-transformation}) the weighting matrix $\mathbf{V}$ is equivariant under the transformations $g$ and $\tilde{g}$, 
\begin{eqnarray}
\mathbf{V}(\tilde{g}(\BS{\beta});\nu^g)
&=&
\int \lambda(\mathbf{f}(\mathbf{z})^\trp\tilde{g}({\BS{\beta}}))^{2} \mathbf{f}(\mathbf{z})\mathbf{f}(\mathbf{z})^\trp\,\nu^g(\mathrm{d} \mathbf{z})
 \nonumber 
 \\
&=&
\int \lambda(\mathbf{f}(g(\mathbf{x}))^\trp\tilde{g}({\BS{\beta}}))^{2} \mathbf{f}(g(\mathbf{x}))\mathbf{f}(g(\mathbf{x}))^\trp\,\nu(\mathrm{d} \mathbf{x})
 \nonumber 
\\
&=&
\int  \lambda(\mathbf{f}(g(\mathbf{x}))^\trp\tilde{g}({\BS{\beta}}))^{2} \mathbf{Q}_g \mathbf{f}(\mathbf{x})\mathbf{f}(\mathbf{x})^\trp \mathbf{Q}_g^\trp\,\nu(\mathrm{d} \mathbf{x})
 \nonumber 
\\
&=&
\mathbf{Q}_g \left(\int \lambda(\mathbf{f}(g(\mathbf{x}))^\trp\tilde{g}({\BS{\beta}}))^{2} \mathbf{f}(\mathbf{x})\mathbf{f}(\mathbf{x})^\trp\,\nu(\mathrm{d}\mathbf{x})\right)  \mathbf{Q}_g^\trp
=
\mathbf{Q}_g \mathbf{V}(\BS{\beta};\nu) \mathbf{Q}_g^\trp ,
\label{eq-transformation-v}
\end{eqnarray}
in the case of a generalized linear model with canonical link. 
This implies
\begin{eqnarray}
\mathrm{IMSE}(\xi^g;\tilde{g}(\BS{\beta}),\nu^g)
&=&
\mathrm{trace}(\mathbf{Q}_g \mathbf{V}(\BS{\beta};\nu)\mathbf{Q}_g^\trp (\mathbf{Q}_g\mathbf{M}(\xi;\BS{\beta})\mathbf{Q}_g^{\trp})^{-1})
\nonumber
\\
&=&
\mathrm{trace}(\mathbf{V}(\BS{\beta};\nu) \mathbf{M}(\xi;\BS{\beta})^{-1})
=\mathrm{IMSE}(\xi;\BS{\beta},\nu) .
\label{eq-transformation-imse}
\end{eqnarray}
Let $\Phi$ be the local \textit{IMSE}-criterion at $\BS{\beta}$ with respect to $\nu$ and $\Phi^{\prime}$ the local \textit{IMSE}-criterion at $\tilde{g}(\BS{\beta)}$ with respect to $\nu^g$, then the \textit{IMSE}-criterion is equivariant under simultaneous transformation of $\BS{\beta}$ and the supplementary argument $\nu$, and by Theorem~\ref{th-equivariance} the locally \textit{IMSE}-optimal design can be transferred.
\begin{corollary}
	\label{cor-equivariance-imse}
	If $\xi^*$ is locally \textit{IMSE}-optimal on $\mathcal{X}$ at $\BS{\beta}$ with respect to $\nu$, then $(\xi^*)^g$  is locally \textit{IMSE}-optimal on $\mathcal{Z}$ at $\tilde{\BS{\beta}}=\tilde{g}(\BS{\beta})$ with respect to $\nu^g$.
\end{corollary}
Note that the results of Corollaries~\ref{cor-equivariance-d} and \ref{cor-equivariance-imse} do not only hold for any generalized linear model, but also, more generally, for all models, where the elemental information matrix is of the form (\ref{eq2-4}) (see e.\,g.\ \citet{schmidtschwabe2017} for further examples).

\begin{exmp*}[Example~\ref{ex-1} continued]
	In order to apply the equivariance result of Corollary~\ref{cor-equivariance-imse} to the gamma model with simple linear regression, $\mathbf{f}(x)=(1,x)^\trp$, the locally \textit{IMSE}-optimal design $\xi^{*}$ on the unit interval $\mathcal{X}=[0,1]$ has to be determined first. 
\end{exmp*}

\begin{proposition} 
	\label{throIMSE} 
	For the one-factor gamma model with simple linear regression $\mathbf{f}(x)^\trp\BS{\beta}=\beta_0+\beta_1 x$ on the experimental region $\mathcal{X}=[0,1]$ locally \textit{IMSE}-optimal designs can be found which are supported on the endpoints $0$ and $1$ of the experimental region.
	\par
	Locally optimal weights $1-w^*$ at $0$ and $w^*$ at $1$, respectively, are given by
	\begin{enumerate}[label=(\alph*)]
		\item $1-w^*=w^*=1/2$ for $\nu$ the uniform (Lebesgue) measure on the interval $[0,1]$,
		\item $1-w^*=(\beta_0+\beta_1)/(2\beta_0+\beta_1)$ and $w^*=\beta_0/(2\beta_0+\beta_1)$ for $\nu$ the (discrete) uniform measure on the endpoints $\{0,1\}$, and
		\item $1-w^*=\beta_0/(2\beta_0+\beta_1)$ and $w^*=(\beta_0+\beta_1)/(2\beta_0+\beta_1)$ for $\nu$ the one-point measure on the midpoint $1/2$ of the design region.
	\end{enumerate}
\end{proposition}

	The proof of Proposition~\ref{throIMSE} is given in the Appendix.
	Note that in Proposition~\ref{throIMSE} the locally optimal weights may depend on the weighting measure $\nu$ used.
	In particular, for the two measures in Proposition~\ref{throIMSE}~(b) and (c) which are concentrated on the endpoints and the midpoint, respectively, the locally optimal weights at $0$ and $1$ are interchanged.
	For the continuous uniform measure (Proposition~\ref{throIMSE}~(a)) equal weights, $w^*=1/2$, are assigned to both endpoints, and the (locally) \textit{IMSE}-optimal design does not depend on the value of the parameter vector $\BS{\beta}$.
	
	Now equivariance can be employed to obtain locally \textit{IMSE}-optimal designs for any other interval $\mathcal{Z}=[a,b]$ as the experimental region.
	We again use the transformation $g(x)=a+cx$, $c=b-a$, together with $\tilde{g}(\BS{\beta})=\mathbf{Q}_g^{-\trp}\BS{\beta}$.
	Let $\xi^*$ be the locally \textit{IMSE}-optimal design of Proposition~\ref{throIMSE} at $\BS{\beta}$ with respect to one of the given weighting measures $\nu$.
	Then, by Corollary~\ref{cor-equivariance-imse}, the design $(\xi^*)^g$ is the locally \textit{IMSE}-optimal design at $\tilde{\BS{\beta}}=\tilde{g}(\BS{\beta})$ with respect to $\nu^g$.
\eop

In order to obtain locally optimal designs at a given value of $\tilde{\BS{\beta}}$ on the transformed design region $\mathcal{Z}$ the inverse transformations $g^{-1}(\mathbf{z})=\mathbf{Q}_g^{-1}\mathbf{z}$ and $\tilde{g}^{-1}(\tilde{\BS{\beta}})=\mathbf{Q}_g\tilde{\BS{\beta}}$ of $g$ and $\tilde{g}$, respectively, have to be used. 
We give this general result only for the case of the \textit{D}- and the \textit{IMSE}-criterion.

\begin{corollary}
	\label{cor-equivariance-reverse}
	Let the equivariance conditions be fulfilled.
	\begin{enumerate}[label=(\alph*)]
		\item
			The design $(\xi^*)^g$ is locally \textit{D}-optimal on $\mathcal{Z}$ at $\tilde{\BS{\beta}}$ if $\xi^*$ is locally \textit{D}-optimal on $\mathcal{X}$ at $\BS{\beta}=\tilde{g}^{-1}(\tilde{\BS{\beta}})$.
		\item
			The design $(\xi^*)^g$ is locally \textit{IMSE}-optimal on $\mathcal{Z}$ at $\tilde{\BS{\beta}}$ with respect to $\nu$ if $\xi^*$ is locally \textit{IMSE}-optimal on $\mathcal{X}$ at $\BS{\beta}=\tilde{g}^{-1}(\tilde{\BS{\beta}})$ with respect to $\nu^{g^{-1}}$.
	\end{enumerate}
\end{corollary}

\begin{exmp*}[Example~\ref{ex-1} continued]
	By Corollary~\ref{cor-equivariance-reverse} we can obtain locally \textit{IMSE}-optimal designs for the one-factor gamma model with simple linear regression $\mathbf{f}(x)^\trp\tilde{\BS{\beta}}=\tilde{\beta}_0+\tilde{\beta}_1 x$ on a given interval $\mathcal{Z}=[a,b]$ with respect to suitably specified weighting measures $\nu_{\mathcal{Z}}$.
	The inversely transformed parameter vector $\BS{\beta}=\tilde{g}^{-1}(\tilde{\BS{\beta}})$ is given by $\BS{\beta}=(\tilde{\beta}_0+a\tilde{\beta}_1,(b-a)\tilde{\beta}_1)^\trp$.
	By Corollary~\ref{cor-equivariance-reverse} and Proposition~\ref{throIMSE} the optimal designs are supported on the endpoints $a$ and $b$ of the interval and the optimal weights $1-w^*$ at $a$ and $w^*$ at $b$, respectively, can be obtained as
	\begin{enumerate}[label=(\alph*)]
		\item 
			$1-w^*=w^*=1/2$ for $\nu_{\mathcal{Z}}$ the uniform (Lebesgue) measure on the interval $[a,b]$,
		\item 
			$1-w^*=(\tilde{\beta}_0+b\tilde{\beta}_1)/(2\tilde{\beta}_0+(a+b)\tilde{\beta}_1)=\mathbf{f}(b)^\trp\tilde{\BS{\beta}}/(\mathbf{f}(a)^\trp\tilde{\BS{\beta}}+\mathbf{f}(b)^\trp\tilde{\BS{\beta}})$ and $w^*=(\tilde{\beta}_0+a\tilde{\beta}_1)/(2\tilde{\beta}_0+(a+b)\tilde{\beta}_1)=\mathbf{f}(a)^\trp\tilde{\BS{\beta}}/(\mathbf{f}(a)^\trp\tilde{\BS{\beta}}+\mathbf{f}(b)^\trp\tilde{\BS{\beta}})$ for $\nu_{\mathcal{Z}}$ the (discrete) uniform measure on the endpoints $\{a,b\}$, and
		\item 
			$1-w^*=(\tilde{\beta}_0+a\tilde{\beta}_1)/(2\tilde{\beta}_0+(a+b)\tilde{\beta}_1)=\mathbf{f}(a)^\trp\tilde{\BS{\beta}}/(\mathbf{f}(a)^\trp\tilde{\BS{\beta}}+\mathbf{f}(b)^\trp\tilde{\BS{\beta}})$ and $w^*=(\tilde{\beta}_0+b\tilde{\beta}_1)/(2\tilde{\beta}_0+(a+b)\tilde{\beta}_1)=\mathbf{f}(b)^\trp\tilde{\BS{\beta}}/(\mathbf{f}(a)^\trp\tilde{\BS{\beta}}+\mathbf{f}(b)^\trp\tilde{\BS{\beta}})$ for $\nu_{\mathcal{Z}}$ the one-point measure on the midpoint $(a+b)/2$ of the experimental region.
	\end{enumerate}
	
	The continuous uniform measure in (a) is the common choice for the \textit{IMSE}-criterion.
	The discrete uniform measure in (b) lays equal interest in the extreme values of the experimental region and may also be applied for the restricted experimental region $\mathcal{X}=\{a,b\}$ which can be used to describe two groups ``$a$'' and ``$b$''.
	In that case the \textit{IMSE}-optimal weights are inverse proportional to the standard deviations $\lambda_x=1/(\mathbf{f}(x)^\trp\tilde{\BS{\beta}})^2$, $x=a,b$, in the groups in accordance with known results on \textit{A}-optimality for group means.
	The one-point measure in (c) coincides with the \textit{c}-criterion for estimating the mean response at the midpoint of the interval.
\eop
\end{exmp*}

Note that the \textit{D}- and \textit{IMSE}-criteria are equivariant with respect to any transformation $g$ of $\mathbf{x}$ for which the regression function $\mathbf{f}$ is linearly equivariant, $\mathbf{f}(g(\mathbf{x}))=\mathbf{Q}_g\mathbf{f}(\mathbf{x})$, and the corresponding transformation $\tilde{g}(\BS{\beta})=\mathbf{Q}_g^{-\trp}\BS{\beta}$ of $\BS{\beta}$.
For other criteria additional requirements may have to be fulfilled by the transformations to obtain equivariance results. 
For example in the case of Kiefer's class of $\Phi_q$-criteria (including the \textit{A-criterion}) the transformation matrix $\mathbf{Q}_g$ should be orthogonal or, at least, satisfy that $\mathbf{Q}_g^{\trp}\mathbf{Q}_g$ is a multiple of the $p\times p$ identity matrix.

For the equivariance of maximin efficiency criteria we require additionally that the underlying local criteria are multiplicatively equivariant with respect to $(g,\tilde{g})$, which means that for every $\BS{\beta}\in\mathcal{B}^{\prime}$ there is a constant $c>0$ such that $\Phi_{\tilde{g}(\BS{\beta})}(\xi^g)=c\Phi_{\BS{\beta}}(\xi)$ uniformly in $\xi$.
Then for the corresponding maximin efficiency criterion we get
\begin{eqnarray}
	\Phi(\xi^g)
	&=&
	\sup_{\tilde{\BS{\beta}}\in\tilde{g}(\mathcal{B}^{\prime})} \frac{\Phi_{\tilde{\BS{\beta}}}(\xi^g)}{\Phi_{\tilde{\BS{\beta}}}(\xi_{\tilde{\BS{\beta}}}^*)}
	\nonumber
	\\
	&=&
	\sup_{\BS{\beta}\in\mathcal{B}^{\prime}} \frac{\Phi_{\tilde{g}(\BS{\beta})}(\xi^g)}{\Phi_{\tilde{g}(\BS{\beta})}((\xi^*_{\BS{\beta}})^g)}
	\nonumber
	\\
	&=&
	\sup_{\BS{\beta}\in\mathcal{B}^{\prime}} \frac{c\Phi_{\BS{\beta}}(\xi)}{c\Phi_{\BS{\beta}}(\xi^*_{\BS{\beta}})}
	= \Phi(\xi) ,
	\label{maximin-eff-equivariance}
\end{eqnarray}
where in the second equality it is used that by Theorem~\ref{th-equivariance} the image of the locally optimal design at $\BS{\beta}$ under $g$ is locally optimal at $\tilde{g}(\BS{\beta})$.
Hence, the resulting maximin efficiency criterion $\Phi$ is equivariant.

By (\ref{eq-transformation-d}) the homogeneous version $\Phi_{\BS{\beta}}(\xi)=(\det(\mathbf{M}(\xi;\BS{\beta})))^{-1/p}$ of the local \textit{D}-criterion is multiplicatively equivariant with $c=\det(\mathbf{Q}_g)^{-2/p}>0$.
Accordingly the local \textit{IMSE}-criterion is multiplicatively equivariant with $c=1$ by (\ref{eq-transformation-imse}).
Hence, both the maximin \textit{D}-efficiency criterion and the maximin \textit{IMSE}-efficiency criterion retain their value under the transformation and are thus equivariant.
\begin{corollary}\label{cor-equivariance-maximin}
	\begin{enumerate}[label=(\alph*)]
		\item
			If $\xi^*$ is maximin \textit{D}-efficient on $\mathcal{B}^{\prime}$, then $(\xi^*)^g$ is maximin \textit{D}-efficient on $\tilde{\mathcal{B}}^{\prime}=\tilde{g}(\mathcal{B}^{\prime})$.
		\item
			If $\xi^*$ is maximin \textit{IMSE}-efficient with respect to $\nu$ on $\mathcal{B}^{\prime}$, then $(\xi^*)^g$ is maximin \textit{IMSE}-efficient with respect to $\nu^g$ on $\tilde{\mathcal{B}}^{\prime}=\tilde{g}(\mathcal{B}^{\prime})$.
	\end{enumerate}
\end{corollary}

\begin{exmp*}[Example~\ref{ex-1} continued]
	In the gamma model with simple linear regression on $[0,1]$ the design $\xi^*$ which assigns equal weights $1/2$ to both endpoints $0$ and $1$ is both locally \textit{D}-optimal and by Proposition~\ref{throIMSE} locally \textit{IMSE}-optimal with respect to the uniform measure $\nu$ on $[0,1]$ for any $\BS{\beta}\in\mathcal{B}$.
	Hence, $\xi^*$ is obviously both maximin \textit{D}-efficient and maximin \textit{IMSE}-efficient with respect to $\nu$ on $\mathcal{B}$ on $[0,1]$.  
	Then, with $g(x)=a+cx$, $c=b-a$, by Corollary \ref{cor-equivariance-maximin} the design $(\xi^*)^g$ which assigns equal weights $1/2$ to $a$ and $b$ is maximin \textit{D}-efficient and maximin \textit{IMSE}-efficient with respect to to the uniform measure $\nu^g$ on  $\tilde{\mathcal{B}}=\tilde{g}(\mathcal{B})$ on $[a,b]$.
\eop
\end{exmp*}

Further maximin \textit{D}- and \textit{IMSE}-efficient designs will be derived in Section~\ref{sec-4}.

\section{Extended equivariance}
\label{sec-3-extended}

As already mentioned in \citet{idais2020analytic} there is a further concept of equivariance in gamma models which is based on the special structure of the intensity function $\lambda$.
First note that for a gamma model with inverse link the maximal region $\mathcal{B}$ of parameter values $\BS{\beta}$ such that the linear component $\mathbf{f}(\mathbf{x})^\trp\BS{\beta}$ is positive constitute a cone in $p$-dimensional Euclidean space, i.\,e.,\ for each vector $\BS{\beta}\in\mathcal{B}$ and every positive scale factor $\tilde{c}>0$ the scaled vector $\tilde{c}\BS{\beta}$ lies also in $\mathcal{B}$.

For given $\tilde{c}>0$, we may consider the (linear) transformation $\tilde{g}(\BS{\beta})=\tilde{c}\BS{\beta}$ of $\BS{\beta}$ in conjunction with the identity mapping $g=\mathrm{id}$ on the experimental region $\mathcal{X}$.
Then, by $\lambda(\mathbf{f}(\mathbf{x})^\trp\tilde{c}\BS{\beta})=\tilde{c}^{\,-2}\lambda(\mathbf{f}(\mathbf{x})^\trp\BS{\beta})$, the intensity is equivariant with respect to $(g,\tilde{g})$ which implies equivariance of the information matrix.

\begin{lemma}
	\label{lemma-2}
	In the gamma model with inverse link the information matrix is equivariant with respect to $g=\mathrm{id}$ and $\tilde{g}(\BS{\beta})=\tilde{c}\BS{\beta}$, i.\,e.,\  
	\[
	\mathbf{M}({\xi;\tilde{c}\BS{\beta}}) = \tilde{c}^{\,-2}\mathbf{M}(\xi;\BS{\beta})
	\] 
	for any positive scale factor $\tilde{c}$.
\end{lemma}

To transfer optimal designs by Theorem~\ref{th-equivariance} it remains to show that the criterion functions under consideration are order preserving with respect to scale transformations.
By Lemma~\ref{lemma-2} we directly get $\det(\mathbf{M}(\xi;\tilde{c}\BS{\beta})) = \tilde{c}^{\,-2p}\det(\mathbf{M}(\xi;\BS{\beta}))$ and, hence, the equivariance of
the \textit{D}-criterion. 
For the \textit{IMSE}-criterion we additionally utilize the equivariance property $\mathbf{V}(\tilde{c}\BS{\beta};\nu) = \tilde{c}^{\,-4}\mathbf{V}(\BS{\beta};\nu)$ of the weighting matrix to get the equivariance, $\mathrm{IMSE}(\xi;\tilde{c}\BS{\beta},\nu) = \tilde{c}^{\,-2}\mathrm{IMSE}(\xi;\BS{\beta},\nu)$.
We thus obtain that designs retain their optimality under scale transformation.

\begin{corollary}
	\label{cor-equivariance-scale}
	In a gamma model on $\mathcal{X}$ with inverse link
	\begin{enumerate}[label=(\alph*)]
		\item
			a locally \textit{D}-optimal design $\xi^*$ at $\BS{\beta}$ is also locally \textit{D}-optimal at $\tilde{\BS{\beta}}=\tilde{c}\BS{\beta}$,
		\item
			a locally \textit{IMSE}-optimal design $\xi^*$ at $\BS{\beta}$ with respect to $\nu$ is also locally \textit{IMSE}-optimal at $\tilde{\BS{\beta}}=\tilde{c}\BS{\beta}$ with respect to $\nu$,
	\end{enumerate}
	for every $\tilde{c}>0$.
\end{corollary}

We may use this result to reduce the number of parameters in the optimization problem.
Therefor we solve the optimization problem first for a standardized parameter setting, where one of the parameters is set to $1$, and then transfer the obtained optimal design to a general parameter vector by setting $\tilde{c}$ to the desired value of the particular parameter used for standardization.

\begin{exmp*}[Example~\ref{ex-1} continued]
	For the one-factor gamma model with simple linear regression $\mathbf{f}(x)^\trp\BS{\beta}=\beta_0+\beta_1 x$ on $\mathcal{X}=[0,1]$ the locally \textit{IMSE}-optimal design at $(1,\gamma)^\trp$ with respect to the discrete uniform weighting measure $\nu$ on the endpoints $\{0,1\}$ assigns weights $1-w^*=1/(2+\gamma)$ and $w^*=(1+\gamma)/(2+\gamma)$ to $0$ and $1$, respectively.
	By Corollary~\ref{cor-equivariance-scale} these weights remain locally optimal for any parameter vector $\BS{\beta}=(\beta_0,\beta_1)^\trp$ with $\beta_1/\beta_0=\gamma$.
	The corresponding reduced parameter region for $\gamma$ is given by $\mathcal{C}=\{\gamma;\, \gamma>-1\}$ which is displayed in Figure~\ref{fig:opt1} as the dashed vertical line at $\beta_1=1$.
	There the diagonal dotted line indicates one ray of values $\BS{\beta}$ which are reduced to one specific value of $\gamma$, i.\,e.,\ $\{\BS{\beta};\, \beta_1=\beta_0=\gamma\}$.
\eop
\end{exmp*}

\begin{exmp*}[Example~\ref{exmp1} continued]
Similarly, in the two-factor gamma model on $[0,1]^2$, the three-dimensional parameter vector $\BS{\beta}=(\beta_0,\beta_1,\beta_2)^\trp$ can be reduced to $\tilde{\BS{\beta}}=(1,\gamma_1,\gamma_2)^\trp$, where $\gamma_1=\beta_1/\beta_0$ and $\gamma_2=\beta_2/\beta_0$, $\beta_0>0$. 
As a consequence, the three-dimensional parameter region $\mathcal{B}$ in Figure~\ref{fig:opt2} is reduced to the two-dimensional region $\mathcal{C}$ for $\BS{\gamma}=(\gamma_1,\gamma_2)^\trp$ in Figure \ref{fig:opt22} which is characterized by the linear constraints $\gamma_1>-1$, $\gamma_2>-1$, and $\gamma_1+\gamma_2>-1$. 
The optimality regions $\mathcal{B}_k$, $k=1,\ldots,4$ of Table \ref{T-1} can now be described in terms of the ratios $\gamma_j=\beta_j/\beta_0$, $j=1,2$, 

\begin{eqnarray}
\mathcal{B}_1 & \equiv & 1- \gamma_1\gamma_2 \leq 0, 
\nonumber
\\
\mathcal{B}_2 & \equiv & (1+\gamma_1+\gamma_2)^2-\gamma_1\gamma_2 \leq 0, 
\nonumber
\\
 \mathcal{B}_3 & \equiv & (1+\gamma_1)^2+\gamma_1\gamma_2 \leq 0, 
\nonumber
\\
 \mathcal{B}_4 & \equiv & (1+\gamma_2)^2+\gamma_1\gamma_2 \leq 0,
\nonumber
\end{eqnarray}
as exhibited in Figure~\ref{fig:opt22}.  
The scaling explains why the subregions in Figure~\ref{fig:opt2} constitute cones in the three-dimensional Euclidean space,
    \begin{figure}[H]
 \centering
    \includegraphics[width=8cm, height=6cm]{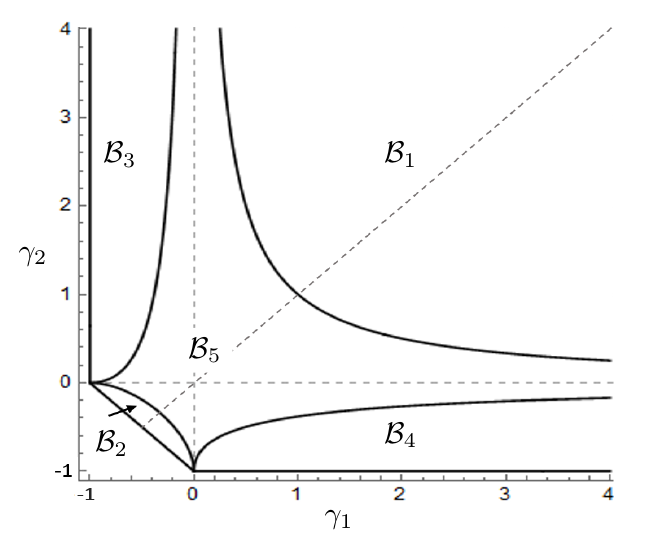}
   \caption{Scaled parameter region of the two-factor gamma model on $[0,1]^2$; the diagonal dashed line represents $\gamma_1=\gamma_2$}
    \label{fig:opt22}
\end{figure}
\eop
\end{exmp*}

By combination of the scale transformation with the transformations of the preceding Section~\ref{sec-3} we get an extended version of Corollaries~\ref{cor-equivariance-d} and \ref{cor-equivariance-imse} by Theorem~\ref{th-equivariance}.
\begin{corollary}
	\label{cor-equivariance-extended}
	In a gamma model on $\mathcal{X}$ with inverse link it holds:
	\begin{enumerate}[label=(\alph*)]
		\item
			If $\xi^*$ is locally \textit{D}-optimal on $\mathcal{X}$ at $\BS{\beta}$, then $(\xi^*)^g$  is locally \textit{D}-optimal on $\mathcal{Z}$ at $\tilde{\BS{\beta}}=\tilde{c}\mathbf{Q}_g^{-\trp}\BS{\beta}$,
		\item
			If $\xi^*$ is locally \textit{IMSE}-optimal on $\mathcal{X}$ at $\BS{\beta}$ with respect to $\nu$, then $(\xi^*)^g$  is locally \textit{IMSE}-optimal on $\mathcal{Z}$ at $\tilde{\BS{\beta}}=\tilde{c}\mathbf{Q}_g^{-\trp}\BS{\beta}$ with respect to $\nu^g$,
	\end{enumerate}
	for every $\tilde{c}>0$.
\end{corollary}
	
This result indicates that for a given transformation $g$ of $\mathbf{x}$ the associated transformation $\tilde{g}$ of $\BS{\beta}$ need not be unique.
Moreover, we may let the scaling factor depend on the parameter vector, $\tilde{c}=\tilde{c}(\BS{\beta})$.
This will lead to a non-linear transformation $\tilde{g}=\tilde{c}(\BS{\beta})\mathbf{Q}_g^{-\trp}\BS{\beta}$ of the parameter vector $\BS{\beta}$.

\begin{exmp*}[Example~\ref{ex-1} continued]
	For the one-factor gamma model with simple linear regression on the unit interval $[0,1]$ consider the reflection $g(x)=1-x$ that maps $[0,1]$ onto itself.
	The corresponding linear transformation of the parameter vector $\BS{\beta}=(\beta_0,\beta_1)^\trp$ is given by $\mathbf{Q}_g^{-\trp}\BS{\beta}=(\beta_0+\beta_1,\,-\beta_1)^\trp$.
	In particular, for a scaled reduced parameter vector $\BS{\beta}=(1,\gamma)^\trp$, $\gamma=\beta_1/\beta_0$, we have $\mathbf{Q}_g^{-\trp}\BS{\beta}=(1+\gamma,\,-\gamma)^\trp$.
	In order to obtain a transformed parameter vector $\tilde{\BS{\beta}}=\tilde{g}(\BS{\beta})$ in reduced form, $\tilde{\beta_0}=1$ the linear transformation has to be rescaled by $\tilde{c}(\BS{\beta})=1/(1+\gamma)$.
	This results in the nonlinear transformation $\tilde{g}((1,\gamma)^\trp)=(1,\,-\gamma/(1+\gamma))^\trp$.
	Note that this transformation is a one-to-one mapping of the maximal region $\mathcal{C}=(-1,\infty)$ for the reduced parameter $\gamma$.
	
	In general, the first entry $\beta_0$ in the parameter vector $\BS{\beta}$ can be preserved by $\tilde{g}(\BS{\beta})=\beta_0/(\beta_0+\beta_1)\mathbf{Q}_g^\trp\BS{\beta}$, where the scaling factor $\tilde{c}(\BS{\beta})=\beta_0/(\beta_0+\beta_1)=1/(1+\gamma)$ only depends on $\gamma=\beta_1/\beta_0$.
\eop
\end{exmp*}

The result of Corollary~\ref{cor-equivariance-extended} carries over directly also for the nonlinear transformation, when $\tilde{c}$ is replaced by $\tilde{c}(\BS{\beta})$.

\begin{exmp*}[Example~\ref{ex-1} continued]
	For the one-factor gamma model with simple linear regression on $[0,1]$ and reflection $g(x)=1-x$ the weighting measures $\nu$ specified in Proposition~\ref{throIMSE} are all invariant with respect to $g$, i.\,e.,\ $\nu^g=\nu$.
	The corresponding locally \textit{IMSE}-optimal designs $\xi^*$ on $[0,1]$ with respect to $\nu$ are supported by the endpoints with optimal weights $1-w^*$ and $w^*$ at $0$ and $1$, respectively.
	Then the designs $(\xi^*)^g$ which assign the exchanged weights $w^*$ to $0$ and $1-w^*$ to $1$ are locally \textit{IMSE}-optimal on $[0,1]$ with respect to $\nu$ at $\tilde{\BS{\beta}}=\tilde{c}(\BS{\beta})\tilde{g}(\BS{\beta})=(\beta_0,\,-\beta_1/(1+\gamma))^\trp$.
\eop
\end{exmp*}

The standardization with respect to the intercept can be extended to more complex models.

\begin{exmp*}[Example~\ref{exmp1} continued]
	In the two-factor gamma model on $[0,1]^2$ we consider the transformation $g_2(\mathbf{x})=(1-x_1,1-x_2)^\trp$ of simultaneous reflection of both explanatory variables and the corresponding rescaled transformation $\tilde{g}_2(\BS{\beta})=\tilde{c}(\BS{\beta})\mathbf{Q}_{g_2}^{-\trp}\BS{\beta}$ of $\BS{\beta}$ which leaves the intercept $\beta_0$ unchanged, i.\,e.,\ $\tilde{c}(\BS{\beta})=\beta_0/(\beta_0+\beta_1+\beta_2)=1/(1+\gamma_1+\gamma_2)$, where $\gamma_j=\beta_j/\beta_0$ for the reduced parameters.
	This induces the transformation $\BS{\gamma}\to -(1/(1+\gamma_1+\gamma_2))\BS{\gamma}$ on the reduced parameter region $\mathcal{C}$ which maps $\mathcal{C}$ onto itself.
	Hence, if a design $\xi^*$ is locally \textit{D}-optimal at $\BS{\gamma}$ which assigns weights $w_i^*$ to $\mathbf{x}_i$, $i=1,\ldots,4$, then the design $(\xi^*)^{g_2}$ which assigns weights $w_4^*$, $w_3^*$, $w_2^*$ and $w_1^*$ to $\mathbf{x}_1,\ldots,\mathbf{x}_4$, respectively, is locally \textit{D}-optimal at $-(1/(1+\gamma_1+\gamma_2))\BS{\gamma}$.
	
	Similar results hold for \textit{IMSE}-optimality.
\eop
\end{exmp*}

For maximin efficiency criteria we additionally allow here that the multiplicative factor in the equivariance of the underlying local criteria
may depend on the parameter $\BS{\beta}$, $c=c_{\BS{\beta}}$.
This does not affect the arguments in (\ref{maximin-eff-equivariance}) and, hence, the resulting maximin efficiency criteria remain equivariant.
The homogeneous version of the local \textit{D}-criterion and the local \textit{IMSE}-criterion are multiplicatively equivariant with $c_{\BS{\beta}}=\tilde{c}(\BS{\beta})^2\det(\mathbf{Q}_g)^{-2/p}>0$ and $c_{\BS{\beta}}=\tilde{c}(\BS{\beta})^{-2}$, respectively.
Hence, for both the maximin \textit{D}-efficiency and the maximin \textit{IMSE}-efficiency criterion their value is not changed under the transformation. 
This criteria are thus equivariant and the result of Corollary \ref{cor-equivariance-maximin} that maximin efficient designs can be transferred remain valid for the nonlinear transformation $\BS{\beta}\to\tilde{c}(\BS{\beta})\mathbf{Q}_g^{-\trp}\BS{\beta}$ in the case of a gamma model.

\section{Invariance} \label{sec-4}

While equivariance can be used to transfer optimal designs the concept of invariance allows to reduce the complexity of finding optimal designs by exploiting symmetries (see e.\,g.\ \citet{schwabe1996optimum}, ch.~3, for the linear model case).
As in linear models we need a (finite) group $G$ of transformation $g$ which map the experimental region $\mathcal{X}$ onto itself.
For each of this transformations $g$ the regression functions $\mathbf{f}$ are assumed to be linearly equivariant, $\mathbf{f}(g(\mathbf{x}))=\mathbf{Q}_g\mathbf{f}(\mathbf{x})$.
For generalized linear models we require additionally that the corresponding transformations of $\BS{\beta}$ also constitute a group such that the pairs of transformations share the group structure.
This requirement is automatically fulfilled for the linear transformations $\tilde{g}(\BS{\beta})=\mathbf{Q}_g^{-\trp}\BS{\beta}$, because the transformation matrices $\mathbf{Q}_g$, $g\in G$, constitute a group with respect to matrix multiplication.
For extended equivariance (Section~\ref{sec-3-extended}) also the scaling factor $\tilde{c}(\BS{\beta})$ has to share the group property which can be seen to hold for the rescaling which leaves the intercept unchanged.

\begin{exmp*}[Example~\ref{ex-1} continued]
	For the one-factor gamma model with simple linear regression on $[0,1]$  the reflection $g(x)=1-x$ maps $[0,1]$ onto itself and is self-inverse, i.\,e.,\ $g^{-1}=g$.
	Hence, $g$ together with the identity $\mathrm{id}$ constitute a group $G=\{\mathrm{id},g\}$ of transformations.
	For $g$ the corresponding transformation of the parameter vector $\BS{\beta}$ is $\tilde{g}(\BS{\beta})=\mathbf{Q}_g^{-\trp}\BS{\beta}=(\beta_0+\beta_1,\,-\beta_1)^\trp$ in the linear case and $\tilde{g}(\BS{\beta}) = \tilde{c}(\BS{\beta})\mathbf{Q}_g^{-\trp}\BS{\beta} = (\beta_0,\,-\beta_0\beta_1/(\beta_0+\beta_1))^\trp$ in the extended case, while in both cases the identity on $\mathcal{B}$ corresponds to the identity $\mathrm{id}$ on $\mathcal{X}$.
	As $\tilde{c}(\tilde{g}(\BS{\beta}))=1+\gamma=1/\tilde{c}(\BS{\beta})$ also $\tilde{g}$ is self-inverse, and the group structure is retained.
\eop
\end{exmp*}

The final ingredient is now that the optimality criterion $\Phi$ is invariant with respect to the group $G$ of transformations, i.\,e.,\ $\Phi(\xi^g)=\Phi(\xi)$ for all $g\in G$ and any design $\xi$.
Then we can use convexity arguments to improve designs by symmetrization.
For this define by $\bar\xi=(1/|G|)\sum_{g\in G}\xi^g$ the symmetrization of a design with respect to the group $G$, where $|G|$ denotes the number of elements in the (finite) group $G$.
Note that $\bar\xi$ is itself a design and it is invariant with respect to $G$, i.\,e.,\ $\bar\xi^g=\bar\xi$ for all $g\in G$. 
If $\Phi$ is invariant and convex we obtain
\begin{equation}
	\Phi(\bar\xi)\leq \frac{1}{|G|}\sum_{g\in G}\Phi(\xi^g)=\Phi(\xi) ,
	\label{eq-convexity}
\end{equation}
where the inequality follows from convexity and the equation from invariance.
From this we can conclude that the invariant designs with respect to $G$ constitute an essentially complete class, which means that we can confine the search for an $\Phi$-optimal design to invariant designs.
\begin{theorem}
	If $\Phi$ is invariant (with respect to $G$) and convex, then there exists an invariant design $\xi^*$ (with respect to $G$) which is $\Phi$-optimal over all designs.
	\label{th-invariance}
\end{theorem}

The class of invariant designs is often much smaller than the class of all designs and optimization an be simplified.
Invariant designs are uniform on orbits $\mathcal{O}_{\mathbf{x}}=\{g(\mathbf{x});\, g\in G\}\subset\mathcal{X}$, i.\,e.,\ all $\mathbf{x}$ in the same orbit have the same weight.
In particular, for an invariant design either all $\mathbf{x}$ in an orbit $\mathcal{O}$ are included with weight $w_{\mathcal{O}}$ or the whole orbit is not in the support of the design.
Then it remains to find the optimal orbits and the corresponding optimal weights which is often a much easier task than to optimize over all possible designs.

\begin{exmp*}[Example~\ref{ex-1} continued]
	For the reflection group $G=\{\mathrm{id},g\}$, $g(x)=1-x$ on $[0,1]$ the orbits are all of the form $\{x,1-x\}$ for $x<1/2$ and $\{1/2\}$ for $x=1/2$, respectively.
	In the one-factor gamma model with simple linear regression it is known that the optimal designs are supported at the endpoints $0$ and $1$ (see \citet{GAFFKE2019}). 
	The only remaining orbit for an optimal design is $\{0,1\}$ and, hence, there is only one invariant design on this orbit which assigns equal weights $1/2$ to each endpoint.
	This design is optimal with respect to each convex invariant criterion.
\eop
\end{exmp*}

In the case of local optimality criteria the requirement of invariance is rather restrictive.
In particular, the local parameter $\BS{\beta}$ has to be invariant under all transformations, i.\,e.,\ $\tilde{g}(\BS{\beta})=\BS{\beta}$ for all $g\in G$.
This condition typically holds only for few values of $\BS{\beta}$.

\begin{exmp*}[Example~\ref{ex-1} continued]
	For the one-factor gamma model with simple linear regression under the reflection $g(x)=1-x$ the parameter $\BS{\beta}$ is only invariant if $\beta_1=0$, i.\,e.,\ there is no effect of the covariate $x$.
	The invariant design which assigns equal weights $1/2$ to the endpoints is locally optimal at $\BS{\beta}$ for $\beta_1=0$.
\eop
\end{exmp*}

\begin{exmp*}[Example~\ref{exmp1} continued]
	For the two-factor gamma model with multiple linear regression on $[0,1]^2$ the reflections $g_2,\ldots,g_4$ are all self-inverse and the composition of each two reflections yield the third.
	Together with the identity $g_1=\mathrm{id}$ the reflections constitute a group $G=\{g_1,g_2,g_3,g_4\}$ of transformations.	
	Locally optimal designs are supported at the vertices of $[0,1]^2$(see \citet{GAFFKE2019}).
	On the vertices there is only one orbit which contains all of them and, hence, the unique invariant design on the vertices assigns equal weights $1/4$ to each of them.
	However, under $G$ the parameter $\BS{\beta}$ is only invariant if $\beta_1=\beta_2=0$, i.\,e.,\ both covariates $x_1$ and $x_2$ do not have an effect.
	Thus the invariant design which assigns equal weights $1/4$ to the vertices is locally optimal at $\BS{\beta}$ only for $\beta_1=\beta_2=0$.
\eop
\end{exmp*}

Note that in both examples above locally optimal designs are obtained for the situation of constant intensity $\lambda$. 
In that case the information matrix is proportional to that in the corresponding linear model with the same linear component and  hence, the locally optimal design coincides with the optimal design for the linear model.

In more complex situations invariance may though be helpful for local optimality at certain parameter values which are invariant with respect to $\tilde{g}$ for all $g\in G$.
To this end first note that in the case of a group $G$ of transformations $g$ the corresponding transformation matrices $\mathbf{Q}_g$ are unimodal, i.\,e.,\ $|\det(\mathbf{Q}_g)|=1$ (see \citet{schwabe1996optimum}, ch.~3).
Hence, by (\ref{eq-transformation-d}) the local \textit{D}-criterion is invariant with respect to $G$ if $\tilde{g}(\BS{\beta})=\BS{\beta}$ for all $g\in G$.
For the \textit{IMSE}-criterion we additionally require that the weighting measure $\nu$ is invariant with respect to $G$, i.\,e.,\ $\nu^g=\nu$ for all $g\in G$.
In that case, by (\ref{eq-transformation-imse}), the local \textit{IMSE}-criterion is invariant with respect to $G$ if $\tilde{g}(\BS{\beta})=\BS{\beta}$ for all $g\in G$.
\begin{corollary}
	\label{cor-invariance-local}
	If $\tilde{g}(\BS{\beta})=\BS{\beta}$ for all $g\in G$, then there exists a locally \textit{D}-optimal design $\xi^*$ at $\BS{\beta}$ which is invariant with respect to $G$.
	
	If additionally $\nu$ is invariant with respect to $G$, then there exists a locally \textit{IMSE}-optimal design $\xi^*$ at $\BS{\beta}$ with respect to $\nu$ which is invariant with respect to $G$.
\end{corollary}

\begin{exmp*}[Example~\ref{exmp1} continued]
	In the two-factor gamma model on $[0,1]^2$ we consider here parameter values $\BS{\beta}$ with $\beta_1=0$, i.\,e.,\ where the first covariate $x_1$ has no effect.
	Such parameter vectors are invariant with respect to the linear transformation $\tilde{g}_3(\BS{\beta})=(\beta_0+\beta_1,-\beta_1,\beta_2)^\trp$ associated with the reflection $g_3(\mathbf{x})=(1-x_1,x_2)^\trp$ of the first covariate $x_1$.
	As the transformation $g_3$ is self-inverse, together with the identity $\mathrm{id}$ it constitutes a group $G_3=\{\mathrm {id}, g_3\}$. 
	Then the local \textit{D}-criterion at such $\BS{\beta}$ with $\beta_1=0$ is invariant with respect to $G_3$, and by Corollary~\ref{cor-invariance-local} a locally \textit{D}-optimal design can be found in the class of designs which are invariant with respect to $G_3$.
	Moreover, also here we can restrict to designs supported by the vertices.
	With respect to $G_3$ the relevant orbits are then $(\mathbf{x}_1,\mathbf{x}_2)$ and $(\mathbf{x}_3,\mathbf{x}_4)$, and invariant designs on the vertices have equal weights $w$ at $\mathbf{x}_1$ and $\mathbf{x}_2$ and weights $1/2-w$ at $\mathbf{x}_3$ and $\mathbf{x}_4$, respectively. 
	We will denote such designs by $\bar\xi_w$.
	The optimization problem for a local \textit{D}-optimal design reduces to finding the optimal weight $w^*$.
	Note that for $\beta_1=0$ the intensities on the orbits are constant, i.\,e.,\ $\lambda_1=\lambda_2$ and $\lambda_3=\lambda_4$, where again $\lambda_i$ denotes the intensity at $\mathbf{x}_i$.
	For the designs $\bar\xi_w$ the determinant of the information matrix becomes 
	\[
		\det(\mathbf{M}(\bar{\xi}_w;\BS{\beta})) 
		= 2(\lambda_1^2\lambda_3 w^2(1/2-w) + \lambda_1\lambda_3^2 w(1/2-w)^2)
	\]
	in the case $\beta_1=0$.
	The optimal weight $w^*$ can be determined by straightforward computations as
	\begin{equation}
		w^* = \frac{3\gamma_2-1+\sqrt{12\gamma_2^2+1}}{6\gamma_2(\gamma_2+2)} 
		\label{opt-we1}
	\end{equation}
	for $\gamma_2=\beta_2/\beta_0 \neq 0$ and $w^*=1/4$ for $\beta_2=0$, and their dependence on  $\gamma_2$ is exhibited in Figure\ref{fig:curve of wbeta1}.
		The resulting invariant design $\bar{\xi}_{w^*}$ is locally \textit{D}-optimal at $\BS{\beta}$ with $\beta_1=0$.
		
		Similar results apply for $\beta_2=0$, when the reflection $g_4$ of the second covariate $x_2$ is used.
\begin{figure}[H]
	\centering
	\includegraphics[width=8cm, height=6.5cm]{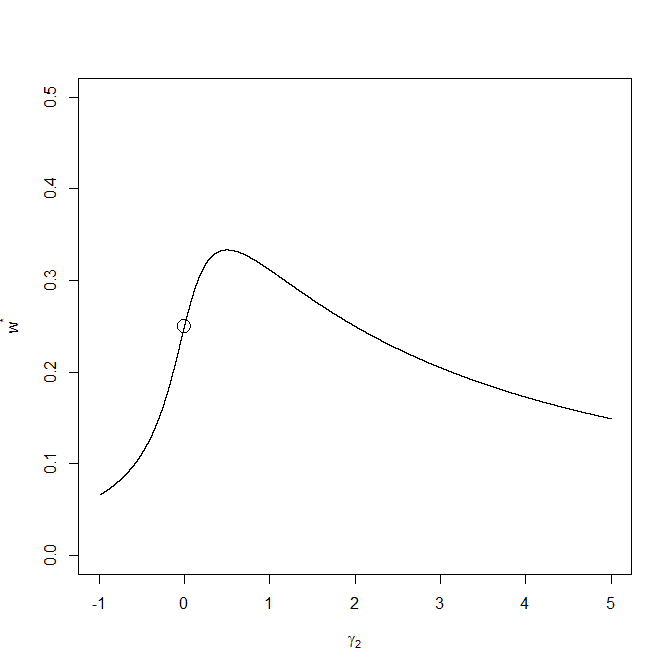}
	\caption{Optimal weights $w^*$ in the two-factor gamma model for $\beta_1=0$}
	\label{fig:curve of wbeta1}
\end{figure}
\eop
\end{exmp*}

\begin{exmp}
	\label{exp-invvv}
	In the two-factor gamma model on $[0,1]^2$ there are also other symmetries which can be employed. 
	In particular, it may be of interest to consider parameter vectors $\BS{\beta}$ with equal slopes,  i.\,e.,\ $\beta_1=\beta_2=\beta$ where the effects of both covariates $x_1$ and $x_2$ are the same.
	These $\BS{\beta}$ are invariant with respect to the linear transformation $\tilde{g}_5(\BS{\beta})=(\beta_0,\beta_2,\beta_1)^\trp$ associated with the permutation $g_5(\mathbf{x})=(x_2,x_1)^\trp$ of the covariates.
	The transformation $g_5$ is self-inverse and constitutes together with the identity $\mathrm{id}$ a group $G=\{\mathrm{id},g_5\}$.
	As locally optimal designs are supported by the vertices of $[0,1]^2$, there are only three relevant orbits $\{\mathbf{x}_1\}$, $\{\mathbf{x}_2,\mathbf{x}_3\}$, and $\{\mathbf{x}_4\}$.
	Optimal invariant designs can thus be characterized by two weight $w_1$ assigned to $\mathbf{x}_1$ and $w_2$ assigned to both $\mathbf{x}_2$ and $\mathbf{x}_3$, while the remaining weight $w_3=1-w_1-2w_2$ is assigned to $\mathbf{x}_4$.
	Note that for $\beta>\beta_0$ and for $-\beta_0/2<\beta \leq \beta_0/3$ minimally supported designs with weights $w_2^*=1/3$ and $w_1^*=1/3$ or $w_3^*=1/3$, respectively, have been seen to be locally \textit{D}-optimal at $\BS{\beta}$ with $\beta_1=\beta_2=\beta$ before.
	In the intermediate case $-\beta_0/3<\beta \leq \beta_0$ the locally \textit{D}-optimal weights can be derived as
	\[
		w_1^*=\frac{3\gamma+1}{4(2\gamma+1)},\  
		w_2^*=\frac{(\gamma+1)^2}{4(2\gamma+1)} \quad\mathrm{ and }\quad
		w_3^*=\frac{1-\gamma}{4},
	\]
	where $\gamma=\beta/\beta_0$  (see \citet{GAFFKE2019}, Theorem~4.3).
	In particular, for $\beta=0$ the uniform weights $w_1^*=w_2^*=w_3^*=1/4$ turn out again to be optimal.

	In the case of the \textit{IMSE}-criterion no explicit formula are available for locally optimal weights so far.
	In Table \ref{T-2} we give some numerical solutions for $\beta_0=1$ and selected values of $\beta_1=\beta_2=\beta$, when $\nu$ is the uniform weighting measure on $[0,1]^2$ which is invariant with respect to $g$.
	These solutions have been obtained by a general non-linear optimization method in the  $\textsf{R}$ package \texttt{Rsolnp} (\citet{R}) using augmented Lagrange multipliers.
	As in the case of the \textit{D}-criterion the locally \textit{IMSE}-optimal design has three support points $\mathbf{x}_1, \mathbf{x}_2,\mathbf{x}_3$, when the effects $\beta$ are large compared to $\beta_0$, but with different weights.
	For intermediate values all four vertices are required.
	In particular, the optimal weights are again uniform on all vertices for $\beta=0$.
	Actually, by the scaling property of Section~\ref{sec-3-extended} the optimal weights depend only on $\gamma=\beta/\beta_0$.

	Moreover, by the additional reflection $g_2(\mathbf{x})=(1-x_1,1-x_2)^\trp$ the optimal weights can be transferred from $\beta>0$ to $\beta<0$ by the nonlinear transformation $\tilde{g}_2$ described in Section~\ref{sec-3-extended}.
	For example, in Table~\ref{T-2} the locally \textit{IMSE}-optimal design at $\tilde{\BS{\beta}}=(1,-3/7,-3/7)^\trp$ is obtained from the locally \textit{IMSE}-optimal design at $\BS{\beta}=(1,3,3)^\trp$ by $g_2$ and the corresponding (nonlinear) transformation $\tilde{g}_2(\BS{\beta})$ in Section~\ref{sec-3-extended}.
\begin{table}[H]
	\caption{Locally \textit{IMSE}-optimal weights in the two-factor gamma model on $[0,1]^2$} 
	\label{T-2}
	\centering
	{
		\scalebox{0.8}{
			\resizebox{\textwidth}{!}{
				\begin{tabular}{rrrllll}
					\hline
					\\[-1ex]
					$\,\,\,\,\,\beta_0$ & $\beta_1$ & $\beta_2$&\,\,\,\,\,\ \  $(0,0)^\trp$  &\,\,\,\,\,\ \ $(0,1)^\trp$ &\,\,\,\,\,\ \ $(1,0)^\trp$&\,\,\,\,\,\ \ $(1,1)^\trp$\\
					\hline
					\\[-2ex]
					\hspace{5ex}1&0&0           &\hspace{3ex}0.250&\hspace{3ex}0.250&\hspace{3ex}0.250&\hspace{3ex}0.250\\
					\hspace{5ex}1&1&1         &\hspace{3ex}0.250&\hspace{3ex}0.300&	\hspace{3ex}0.300 &\hspace{3ex}0.150    \\
					\hspace{5ex}1&2&2         &\hspace{3ex}0.242&\hspace{3ex}0.362&	\hspace{3ex}0.362 &\hspace{3ex}0.034         \\
					\hspace{5ex}1&3&3         &\hspace{3ex}0.236&\hspace{3ex}0.382&	\hspace{3ex}0.382 &\hspace{3ex}0.000    \\
					\hspace{5ex}1&10&10        &\hspace{3ex}0.214&\hspace{3ex}0.393&	\hspace{3ex}0.393 &\hspace{3ex}0.000     \\
					\hspace{5ex}1&$-3/7$&$-3/7$        &\hspace{3ex}0.000&\hspace{3ex}0.382&	\hspace{3ex}0.382 &\hspace{3ex}0.236       \\
					\hline
		\end{tabular}}}
	}
\end{table}
\eop
\end{exmp}

We now turn to maximin efficiency where invariance can become a powerful tool.
For this we additionally require that the subregion $\mathcal{B}^{\prime}$ of interest is also invariant with respect to the pairs $(g,\tilde{g})$ of transformations, i.\,e.,\ $\tilde{g}(\mathcal{B}^{\prime})=\mathcal{B}^{\prime}$ for all $g\in G$.
As already mentioned in Section~\ref{sec-3} the \textit{D}- and the \textit{IMSE}-criterion are multiplicatively equivariant and, hence, by (\ref{maximin-eff-equivariance}) their value is not changed under the transformations. 
Thus both maximin \textit{D}- and \textit{IMSE}-efficiency are invariant with respect to any group $G$ of transformations.
\begin{corollary}
	\label{cor-invariance-maximin}
	If $\tilde{g}(\mathcal{B}^{\prime})=\mathcal{B}^{\prime}$ for all $g\in G$, then there exists a maximin \textit{D}-efficient design $\xi^*$ on $\mathcal{B}^{\prime}$ which is invariant with respect to $G$.
	
	If additionally $\nu$ is invariant with respect to $G$, then there exists a maximin \textit{IMSE}-efficient design $\xi^*$ on $\mathcal{B}^{\prime}$ with respect to $\nu$ which is invariant with respect to $G$.
\end{corollary}

\begin{exmp*}[Example~\ref{ex-1} continued]
	In the one-factor gamma model with simple linear regression on $[0,1]$ the invariant design $\bar\xi$ which assigns equal weights $1/2$ to the endpoints is both maximin \textit{D}- and \textit{IMSE}-efficient with respect to the uniform measure $\nu$ on $[0,1]$ on $\mathcal{B}$ as has already been pointed out at the end of Section~\ref{sec-3}.
	
	In general, however, in contrast to the local criteria there is no direct majorization argument available for maximin efficiency criteria which allows to restrict the support of an optimal design to the extremal points of the experimental region.
	Therefore, to keep argumentation simple and to concentrate on the concept of invariance, we deliberately confine the support to these endpoints.
	Then with respect to reflection, $g(x)=1-x$, the only invariant design which assigns equal weights $1/2$ to the endpoints, is maximin efficient for any invariant criterion.
	In particular, this design is maximin \textit{IMSE}-efficient on $\mathcal{B}$ with respect to any invariant weighting measure $\nu$ as specified in Proposition~\ref{throIMSE}.
\eop
\end{exmp*}

Although, in general, the value of the maximin efficiency criterion requires the knowledge of the locally optimal designs, the maximin efficient design may be constructed without this information as the above example shows.
The results can be extended to more complex models.

\begin{exmp*}[Example~\ref{exmp1} continued]
	In the two-factor gamma model on $[0,1]^2$ with multiple regression we first consider maximin efficiency on the region $\mathcal{B}$ of all possible values for the parameter vector.
	This region is invariant under the transformations associated with the group $G=\{g_1,\ldots,g_4\}$ of reflections of the covariates.
	Similar to the case of  the one-factor gamma model we deliberately confine the support of the designs to the vertices $\mathbf{x}_1,\ldots,\mathbf{x}_4$ of the experimental region to keep argumentation simple.
	Then there is only one orbit which contains all vertices, and the only invariant design with respect to $G$ is the uniform design $\bar\xi$ on the vertices which assigns equal weights $1/4$ to each vertex.
	Hence, the design $\bar\xi$ is maximin efficient on $\mathcal{B}$ for any invariant criterion with respect to $G$.
	
	This result carries over to any parameter subregion $\mathcal{B}^{\prime}$ which is invariant with respect to $G$.
	For example, if the intercept $\beta_0$ is restricted to a subset, $\beta_0\in\mathcal{B}_0$, of its marginal region $(0,\infty)$ or set to a fixed value ($\mathcal{B}_0=\{\beta_0\}$), while the slopes may vary across their corresponding (conditional) regions, then the resulting subregion $\mathcal{B}^{\prime}=\{\BS{\beta}\in\mathcal{B};\, \beta_0\in\mathcal{B}_0\}$ is invariant with respect to the rescaled transformation $\tilde{g}$ associated with $G$.
	Hence, the uniform design $\bar\xi$ is also maximin efficient on $\mathcal{B}^{\prime}$ for any invariant criterion with respect to $G$.
	In particular, this holds for the reduced parameter region $\mathcal{C}$ displayed in Figure~\ref{fig:opt22} for $\beta_0=1$ fixed.
\eop
\end{exmp*}

Invariance can also be employed in cases, where there are less symmetries and thus there are more than one orbit such that still the weights between the orbits have to be optimized.

\begin{exmp*}[Example~\ref{exp-invvv} continued]
	In the two-factor gamma model on $[0,1]^2$ now main interest is in the parameter vectors $\BS{\beta}$ with equal slopes, i.\,e.,\ $\beta_1=\beta_2=\beta$ we consider the parameter subregion $\mathcal{B}^{\prime}=\{(\beta_0,\beta,\beta)^\trp;\, \beta>-\beta_0/2,\beta_0>0\}$.
	In terms of the reduced parameter $\BS{\gamma}=(\gamma_1,\gamma_2)^\trp$, $\gamma_j=\beta_j/\beta_0$, the subset $\mathcal{B}^{\prime}$ reduces to $\mathcal{C}^{\prime}=\{(\gamma,\gamma)^\trp: \gamma>-1/2\}$ which is exhibited in Figure~\ref{fig:opt22} by the dashed diagonal line. 

	For the transformation $g_2(\mathbf{x})=(1-x_1,1-x_2)^\trp$ of simultaneous reflection of both explanatory variables we now consider the rescaled transformation $\tilde{g}_2(\BS{\beta})=\tilde{c}(\BS{\beta})\mathbf{Q}_{g_2}^{-\trp}\BS{\beta}$ of the parameter vector $\BS{\beta}$ which leaves the intercept $\beta_0$ unchanged, i.\,e.,\ $\tilde{c}(\BS{\beta})=\beta_0/(\beta_0+\beta_1+\beta_2)=1/(\gamma_0+\gamma_1+\gamma_2)$ which becomes $\tilde{c}(\BS{\beta})=\beta_0/(\beta_0+2\beta)=1/(1+2\gamma)$ on the subsets $\mathcal{B}^{\prime}$ and $\mathcal{C}^{\prime}$, respectively.
	In particular, the relevant reduced slope parameter $\gamma>-1/2$ is mapped to $-\gamma/(1+2\gamma)>-1/2$.
	Obviously, both $\mathcal{B}^{\prime}$ and $\mathcal{C}^{\prime}$ are invariant with respect to $\tilde{g}_2$.
	
	To make use of the symmetries with respect to the transformations $g_2$ and $g_5$ together we consider the group $G^{\prime}=\{\mathrm{id},g_2,g_5,g_6\}$ generated by them, where the composition $g_6$ of $g_2$ and $g_5$ is the reflection at the secondary diagonal of the unit square, $g_6(\mathbf{x})=(1-x_2,1-x_1)^\trp$.
	We restrict again to the vertices of the experimental region.
	Then there are two orbits  $\{\mathbf{x}_1,\mathbf{x}_4\}$ and $\{\mathbf{x}_2,\mathbf{x}_3\}$. 
	Optimal invariant designs $\bar\xi_w$ can thus be characterized by a single quantity $w$, $0<w<1/2$, where $w_1=w$ is the weight assigned to each of $\mathbf{x}_1$ and $\mathbf{x}_4$ in the first orbit and $w_2=1/2-w$ is the weight assigned to each of $\mathbf{x}_2$ and $\mathbf{x}_3$ in the second orbit.
	Optimization is then reduced to determine the optimal weight $w$.
	
	For the \textit{D}-criterion the determinant of an invariant design $\bar{\xi}_w$ is given by
	\[
		\det (\mathbf{M}(\bar{\xi}_w;\BS{\beta})) 
		= \frac{w(1-2w)\left((1+\gamma)^2+\gamma^2(1-2w)\right)}{2\beta_0^6(1+\gamma)^4(1+2\gamma)^2}
	\]
	locally at $\BS{\beta}=(\beta_0,\beta,\beta)^\trp$, where $\gamma=\beta/\beta_0$.
	
	To find the maximin \textit{D}-efficient design over $\mathcal{B}^{\prime}$ we can confine the analysis to the reduced parameter region $\mathcal{C}^{\prime}$, i.\,e.,\ $\gamma>-1/2$.
	For $\gamma\geq 1$ the minimally supported design $\xi^*_{\BS{\beta}}$ with equal weights $1/3$ on $\mathbf{x}_1, \mathbf{x}_2, \mathbf{x}_3$ is locally \textit{D}-optimal and has $\det(\mathbf{M}(\xi_{\BS{\beta}}^*;\BS{\beta}))=\beta_0^6(1+\gamma)^4/27$.
	Hence for the \textit{D}-efficiency of $\bar{\xi}_w$ at $\gamma\geq 1$ we get
	\begin{equation}
		\mathrm{eff}_{\textit{D}}(\bar{\xi}_w;\BS{\beta})^3 = 27\frac{w(1-2w)\left((1+\gamma)^2+\gamma^2(1-2w)\right)}{2(1+2\gamma)^2} .
		\label{eff-fun}
	\end{equation}
	This expression is decreasing in $\gamma\geq 1$. 
	Therefore, its infimum is obtained when $\gamma$ tends to $\infty$ and the limiting values is $h(w)=\inf_{\gamma\geq 1} \mathrm{eff}_{\textit{D}}(\bar{\xi}_w;\BS{\beta})^3=27w(1-2w)(1-w)/4$.
	The function attains it maximum over $0<w<1/2$ at $w^*=(3-\sqrt{3})/6=0.2113$. 
	The minimal efficiency of $\bar{\xi}_{w^*}$ over $\gamma\geq 1$ is then $\mathrm{eff}_{\textit{D}}(\bar{\xi}_w;\BS{\beta})=h(w^*)^{1/3}=0.8660$, and $\bar{\xi}_{w^*}$ is maximin \textit{D}-efficient over $\gamma\geq 1$ within the class of invariant designs $\bar{\xi}_w$.  
	By symmetry considerations with respect to the transformation $g_2$ (or $g_6$) this result carries over to the region $-1/2<\gamma\leq -1/3$.
	For the intermediate region ($-1/3<\gamma<1$) the \textit{D}-efficiencies of $\bar\xi_{w^*}$ have been computed numerically and are displayed in Figure \ref{fig:D-eff}.
	As the \textit{D}-efficiency of $\bar\xi_{w^*}$ is increasing in $\gamma$ for $\gamma<0$ and decreasing for $\gamma>0$ its minimal efficiency is attained at the boundary of the parameter region and, hence, it can be concluded that $\bar\xi_{w^*}$ is maximin \textit{D}-efficient on $\mathcal{C}^{\prime}$ with minimal \textit{D}-efficiency $0.8660$.
	This result carries over to the whole region $\mathcal{B}^{\prime}$ of equal slopes as well as to subregions $\{\BS{\beta}\in\mathcal{B}^{\prime};\, \beta_0\in\mathcal{B}_0\}$ with constraints on the intercept.
	For comparison, in Figure~\ref{fig:D-eff} the \textit{D}-efficiency is also plotted for the uniform design $\bar{\xi}_{1/4}$ which is locally optimal at $\gamma=0$.
	By (\ref{eff-fun}) the maximin \textit{D}-efficiency of the uniform design $\bar{\xi}_{1/4}$ can be computed as $0.8585$ which is slightly worse than the value for the maximin \textit{D}-efficient design $\bar\xi_{w^*}$.
 \begin{figure}
 \centering
    \includegraphics[width=8cm, height=6.5cm]{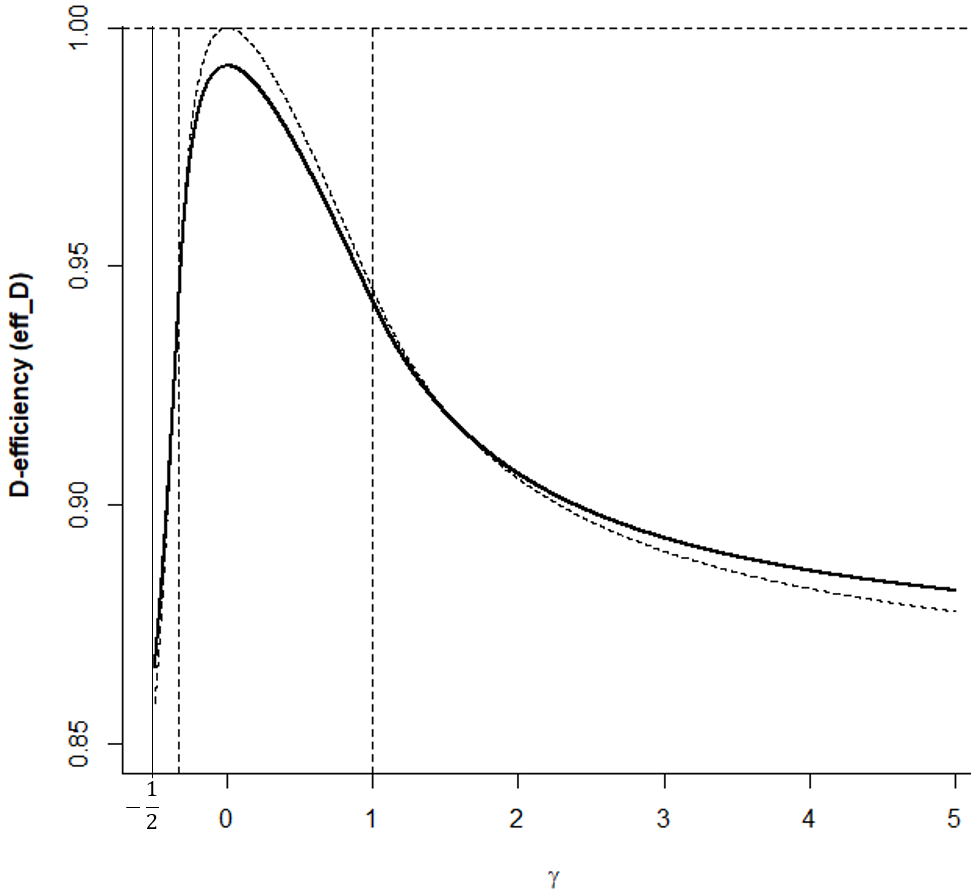}
   \caption{\textit{D}-efficiency of the maximin \textit{D}-efficient design $\bar \xi_{w^*}$ (solid line) and the uniform design $\bar \xi_{1/4}$ (dashed line); vertical lines indicate the lower bound ($\gamma=-1/2$) and the thresholds $\gamma=-1/3$ and $\gamma=1$}
    \label{fig:D-eff}
\end{figure}
\eop
\end{exmp*}


\section{Conclusion and discussion}
\label{sec-5}
In this article we present an overview of the concept of in- and equivariance in the design of experiments for generalized linear models. 
In contrast to the well-known results in linear models, where only the experimental settings are transformed, we have to consider pairs of transformations in generalized linear models which act simultaneously on the experimental settings and on the location parameters in the linear component.
We focus on local optimality and maximin efficiency for the common \textit{D}- and \textit{IMSE}-criteria which allow a wide range of transformations for the experimental settings like scaling, permutations or reflections.
For other criteria the use of invariance is rather limited, because additional structures of the transformations would be required like orthogonality.

As in linear models the transformation of the experimental settings has to act in a linear way on the regression functions of the linear component.
The parameters can then be transformed linearly in such a way that the value of the linear component and, hence, of the intensity is not changed (see \citet{radloff2016invariance}).
Besides this natural choice also nonlinear transformations of the parameters may be employed, if additional properties of the intensity function can be used. 
We illustrate this feature by the gamma model with inverse link for which the intensity is only scaled by a multiplicative factor based on the parameter.
This scaling does not affect standardized design criteria like maximin efficiency, and invariance can also be used here. 
All results can be extended directly to other model specifications like censoring when the intensity depends only on the linear component (see \citet{schmidtschwabe2017} for examples).

For the case of maximin efficiency it would be desirable to obtain also there majorization results as for local optimality which allow for reducing the experimental settings to the extremal points of the experimental region.
However, the findings in \citet{GAFFKE2019} do not carry over, because the arguments used there are of a local nature and do not work uniformly on a parameter region.
Alternatively, equivalence theorems could be employed for maximin efficiency (see \citet{pronzato-pazman}, ch.~8), but they require that the minimal efficiency is attained in the parameter region which is violated in the example.
Therefore the restriction to the extremal points of the experimental region remains an open problem.

\section*{Acknowledgment}
The first author acknowledges support by the the emergency fund of the Graduate Academy, Otto-von-Guericke-University Magdeburg, on the occasion of the COVID-19 pandemic.

\section*{Confict of interest}
On behalf of all authors, the corresponding author states that there is no conflict of interest.

 \section*{Appendix}
\hspace{-4ex} {\textbf{Proof of Proposition~\ref{throIMSE} }: } 
First note that by the majorization argument in \citet{GAFFKE2019} there is an optimal design of the form $\xi_w$ which assigns weight $w$ to $1$ and weight $1-w$ to $0$.
The information matrix for $\xi_w$ and its inverse are
\[
	\mathbf{M}(\xi_{w};\BS{\beta})
		= \left(\begin{array}{cc}
			(1-w)/a^2+w/b^2&w/b^2
			\\
			w/b^2&w/b^2
		\end{array}\right)
	\quad \mathrm{ and } \quad
	\mathbf{M}(\xi_{w};\BS{\beta})^{-1}
		= \left(\begin{array}{cc}
			a^2/(1-w)&-a^2/(1-w)
			\\
			-a^2/(1-w)&a^2/(1-w)+b^2/w
		\end{array}\right) ,
\]
where $a=\beta_0$ and $b=\beta_0+\beta_1$ for simplification. 

For given measure $\nu$ denote by $v_k=\int x^k/(\beta_0+\beta_1 x)^{4}\nu({\rm d} x)$, $k=0,1,2$ the entries in the weighting matrix $\mathbf{V}(\BS{\beta};\nu)=\left(\begin{array}{cc}v_0&v_1\\v_1&v_2\end{array}\right)$.
Then the \textit{IMSE}-criterion becomes
\[
	\mathrm{IMSE}(\xi_w;\BS{\beta},\nu)
	=\mathrm{trace}(\mathbf{V}(\BS{\beta};\nu)\mathbf{M}(\xi_{w};\BS{\beta})^{-1})
	=(1/(1-w))a^2(v_0-2v_1+v_2)+(1/w)b^2 v_2 .
\]
To obtain the $v_k$ we introduce the auxiliary terms $c_k=\int 1/(\beta_0+\beta_1 x)^{k}\nu({\rm d} x)$, $k=2,3,4$.
Then the $v_k$ can be represented as
\[
v_0=c_4, 
\quad v_1=(c_3-a c_4)/(b-a) 
\quad \mathrm{ and } \quad v_2=(c_2-2a c_3+a^2 c_4)/(b-a)^2 ,
\]
 if $\beta_1=b-a\neq 0$.
 The case $\beta_1=0$, when the intensity is constant, follows directly from the corresponding linear model.
\begin{enumerate}[label=(\alph*)]
\item 
	For the (continuous) uniform measure $\nu$ on $[0,1]$ we have $c_k=(a^{-k+1}-b^{-k+1})/((k-1)(b-a))$ and, hence, 
	\[
		v_0=\frac{a^2+ab+b^2}{3a^3b^3}, 
		\quad v_1=\frac{2a+b}{6a^2b^3} 
		\quad \mathrm{ and } \quad v_2=\frac{1}{3ab^3} .
	\]
	From this we get 
	\[
		\mathrm{IMSE}(\xi_w;\BS{\beta},\nu)
		=\frac{1}{3ab}\left(\frac{1}{w}+\frac{1}{1-w}\right) 
	\]
	which is optimized for $w^*=1/2$.
\item  
	For the (discrete) uniform measure $\nu$ on $\{0,1\}$ we have
	$c_k=(a^{-k}+b^{-k})/2$ and, hence, 
	\[
	v_0=(a^{-4}+b^{-4})/2, 
	\quad \mathrm{ and } \quad v_1=v_2=b^{-4}/2 .
	\]
	From this we get 
	\[
		\mathrm{IMSE}(\xi_w;\BS{\beta},\nu)
		=\frac{1}{2}\left(\frac{1}{b^2 w}+\frac{1}{a^2(1-w)}\right) 
	\]
	which is optimized for $w^*=a/(a+b)$.
\item 
	For the one-point measure $\nu$ at $1/2$ we have
	$c_k=\frac{2^k}{(a+b)^k}$ and, hence, $v_k=2^{4-k}/(a+b)^4$.
	From this we get 
	\[
		\mathrm{IMSE}(\xi_w;\BS{\beta},\nu)
		=\frac{4}{(a+b)^4}\left(\frac{a^2}{1-w}+\frac{b^2}{w}\right) 
	\]
	which is optimized for $w^*=b/(a+b)$.
\end{enumerate}
\end{document}